\documentclass[11pt]{amsart}
\usepackage[margin=1in]{geometry}
\usepackage{amsmath,amssymb,amsthm,mathtools}
\usepackage{booktabs,array,multirow}
\usepackage{enumitem}
\usepackage{microtype}
\usepackage[hidelinks]{hyperref}
\usepackage{xurl, graphics,xcolor}
\definecolor{additionblue}{RGB}{20,100,165}
\hypersetup{pdftitle={Foliations, free cores, and Terao's conjecture for seventeen lines},pdfauthor={Felipe Monteiro and Henrique N. Sa Earp and Tomas S. R. Silva},pdfsubject={Terao's conjecture for seventeen complex projective lines}}
\newcommand{\ghlogo}{\raisebox{-0.2em}{\includegraphics[height=1em]{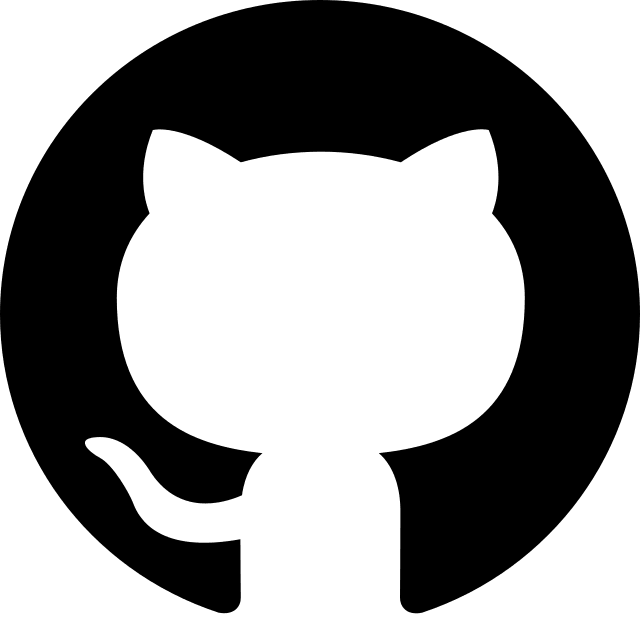}}}
\newtheorem{theorem}{Theorem}[section]
\newtheorem{lemma}[theorem]{Lemma}
\newtheorem{proposition}[theorem]{Proposition}
\newtheorem{corollary}[theorem]{Corollary}
\theoremstyle{definition}
\newtheorem{definition}[theorem]{Definition}
\newtheorem{example}[theorem]{Example}
\theoremstyle{remark}
\newtheorem{remark}[theorem]{Remark}
\newcommand{\C}{\mathbb C}
\newcommand{\Q}{\mathbb Q}

\newcommand{\PP}{\mathbb P}
\newcommand{\cA}{\mathcal A}

\newcommand{\cE}{\mathcal E}
\newcommand{\cF}{\mathcal F}
\newcommand{\cO}{\mathcal O}
\newcommand{\cI}{\mathcal I}
\newcommand{\chib}{\bar\chi}
\newcommand{\eu}{E}
\newcommand{\code}[1]{\texttt{#1}}
\DeclareMathOperator{\Der}{Der}
\DeclareMathOperator{\mdr}{mdr}
\DeclareMathOperator{\Bour}{Bour}
\DeclareMathOperator{\Sing}{Sing}
\DeclareMathOperator{\mult}{mult}
\DeclareMathOperator{\ord}{ord}

\setlist[enumerate]{leftmargin=2em,itemsep=2pt,topsep=4pt}
\setlist[itemize]{leftmargin=1.7em,itemsep=2pt,topsep=4pt}

\title[Foliations, free cores, and Terao's conjecture for seventeen lines]{Foliations, free cores, and Terao's conjecture for seventeen lines}
\author{Felipe Monteiro}
\author{Henrique N. S\'a Earp}
\author{Tom\'as S. R. Silva}
\address[F. Monteiro, H. N. S\'a Earp, T. S. R. Silva]{Instituto de Matem\'atica, Estat\'istica e Computa\c{c}\~ao Cient\'ifica (IMECC), Universidade Estadual de Campinas (UNICAMP), 13083-859, Campinas, S\~ao Paulo, Brazil}
\address[F. Monteiro]{Universite Bourgogne Europe, CNRS, IMB UMR 5584, F-21000 Dijon, France}
\email[F. Monteiro]{ffcmonteiro00@gmail.com}
\email[H. N. S\'a Earp]{henrique.saearp@ime.unicamp.br}
\email[T. S. R. Silva]{tomas@ime.unicamp.br}

\date{\today}
\subjclass[2020]{Primary 14N20; Secondary 32S22, 52C35, 05B35}
\keywords{Line arrangement, Terao's conjecture, logarithmic bundle, foliation, nearly free arrangement, computer-assisted proof}
\begin{document}
\begin{abstract}
We prove that freeness is determined by the intersection lattice for arrangements of seventeen lines in the complex projective plane. The proof reduces a possible counterexample to smaller arrangements whose realisation spaces can be computed. A minimal logarithmic derivation defines a foliation; its inflection divisor and a local singularity count imply that the non-free member of a counterexample is nearly free with exponents $(7,10)$. Deleting four lines then produces a thirteen-line arrangement whose second Chern class is controlled by a combinatorial deficit. When the deficit is four, the remaining arrangement is free with exponents $(6,6)$. An exact covering argument places the necessary thirteen-line cores in four finite classes, containing $3{,}829{,}691$ lattices. We compute their realisation spaces and reconstruct all relevant seventeen-line extensions, including specialisations of positive-dimensional families. Among the candidate lattices, exactly sixteen admit a non-free realisation; they are rigid and have no free realisation. Together with the known cases, this establishes Terao's conjecture for central arrangements of at most seventeen planes in $\C^3$.\\
\small Code: \ghlogo \, \url{https://github.com/TomasSilva/TC17}~~\cite{TC17}.
\end{abstract}
\maketitle

\section{Introduction}\label{sec:intro}

An arrangement of lines carries two kinds of information. Its intersection lattice records which lines meet at which points. Its module of logarithmic derivations depends on the equations of the lines. Terao's conjecture asks whether the freeness of that module can nevertheless be read from the lattice. We work over $\C$ and prove the following case.

\begin{theorem}\label{thm:main}
If two arrangements of seventeen lines in $\PP^2_\C$ have isomorphic intersection lattices, then either both are free or neither is free. Consequently, freeness of a central arrangement of at most seventeen planes in $\C^3$ is determined by its intersection lattice.
\end{theorem}

The cases of at most twelve, thirteen, and fourteen lines were established in \cite{FV14,DIM19,BBJKL21}; Barakat and K\"uhne subsequently treated at most fourteen hyperplanes over every field \cite{BK23}. Anghel and Stan proved the complex fifteen- and sixteen-line cases in \cite{AS26}. We obtained those two cases independently in computations recorded on 5--6 October 2026 in \cite{TeraoHunter}, before learning of their work on 20 October. The two computations agree on their common range. We give no separate account of those cases here and use \cite{AS26} for the theorem up to sixteen lines.

The difficulty at seventeen lines is not deciding freeness for a given arrangement. It is finding a finite collection of arrangements that contains every possible counterexample. Direct enumeration under the usual numerical constraints becomes too large: sampling the candidate classes gives an estimated order of $10^{11}$ lattices. Our reduction instead asks which smaller arrangements must occur inside a counterexample and how the missing lines can be recovered, reducing the needed computer verification to just $3{,}829{,}691$ lattices of thirteen lines.

Suppose that $\cA_1$ is free, $\cA_2$ is not, and their seventeen-line lattices agree. Three observations make the smaller-arrangement approach effective. First, the least-degree Jacobian syzygy of $\cA_2$ gives a foliation preserving all its lines. The inflection divisor has degree $3r$, where $r$ is the syzygy degree. A count of the singularities shows that its residual divisor has degree at least two, see Theorem~\ref{thm:foliation}. Thus $3r\ge19$, which forces the free exponents of $\cA_1$ to be $(1,8,8)$ and $\cA_2$ to be nearly free with exponents $(7,10)$, by Theorem~\ref{thm:17}.

Second, there is a useful choice of four lines to delete. For a set $S$ of four lines, let $s_H$ count the intersection points on $H$ and let $f(S)$ count the points whose incident lines all belong to $S$. The quantity
\[
\varphi(S)=\sum_{H\in S}(s_H-7)-f(S)
\]
measures the change in the second Chern class of the logarithmic bundle. It is at most four. When it equals four, the remaining thirteen lines form a free arrangement with exponents $(6,6)$; this follows from a short Riemann--Roch argument. If an eight-point line exists, a finite covering argument guarantees $\varphi(S)\ge1$ for some $S$. The thirteen-line cores therefore belong to four enumerable classes. The complementary case, in which every line has at most seven points, is small enough to handle directly.

Third, the deletion loses no need to consider a continuous choice of an additional line: each missing line must join two intersection points of the arrangement already present, by Lemma~\ref{lem:determined}. We can therefore extend the realisation spaces of the cores exactly. This is the main computational distinction from the previous proofs for fifteen and sixteen lines: the search runs through smaller realisation spaces and their extensions, rather than through all candidate seventeen-line lattices.

Sections~\ref{sec:prelim}--\ref{sec:seventeen} establish the numerical and geometric reduction. Section~\ref{sec:RR} explains the deletion argument, Section~\ref{sec:cover} proves the covering statement, and Section~\ref{sec:computation} completes the finite verification. The computation and its checks are documented in Appendix~\ref{sec:data} and in \cite{TC17}. In related independent work, Anghel and Stan use foliations to prove combinatorial invariance of the generic splitting type for at most sixteen lines \cite{AS26split}; the argument here is directed at the nearly free seventeen-line case and its thirteen-line cores.

\subsection*{Acknowledgements} We thank Cristian Anghel and Rare\c{s} Stan for their helpful correspondence and valuable discussions. \textbf{FM} is supported by São Paulo Research Foundation (FAPESP) grant 2021/10550-4. \textbf{HNSE} is supported by FAPESP grants 2020/09838-0, with BI0S--Brazilian Institute of Data Science;
2021/04065-6, with BRIDGES--Brazil--France Interplays in Gauge Theory, Extremal Structures and
Stability; 2024/00923-6, with CBG--Brazilian Centre for Geometry; and by the Brazilian National
Council for Scientific and Technological Development (CNPq), through the PQ-A Productivity Grant
307145/2025-5. \textbf{TSRS} is supported by FAPESP grant 2022/09891-4. This work was supported in part by a grant of access to OpenAI models through the ChatGPT for Academic Researchers program.

\subsection*{Data Availability}
The code and datasets used in this work are available at: 
\begin{center}
\url{https://github.com/TomasSilva/TC17}~~\cite{TC17}.
\end{center}

\subsection*{AI Usage Statement} Claude (Anthropic) assisted with the computations and the first draft of this article. 
Codex (OpenAI) assisted with revising the exposition, organisation, internal references, bibliography and \LaTeX{} source, and with checking mathematical arguments against the cited sources, selected numerical examples, compilation and layout. These revisions were shepherded and validated by the authors, who take responsibility for the content. The codebase was also developed with assistance from Claude (Anthropic).

\section{Logarithmic derivations and numerical constraints}\label{sec:prelim}

\subsection{The bundle and its Chern classes}
Let $R=\C[x,y,z]$ and let $Q=\prod_{H\in\cA}\alpha_H$ define an arrangement of $n$ distinct lines. A polynomial vector field is logarithmic along $\cA$ if it is tangent to each line:
\[
D(\cA)=\{\theta\in\Der_\C(R):\alpha_H\mid\theta(\alpha_H)\text{ for every }H\in\cA\}.
\]
Writing $\eu=x\partial_x+y\partial_y+z\partial_z$, the identity $\eu(Q)=nQ$ gives
\[
D(\cA)=R\eu\oplus D_0(\cA),\qquad
D_0(\cA)=\{(a,b,c)\in R^3:aQ_x+bQ_y+cQ_z=0\}.
\]
Indeed, subtracting $(\theta(Q)/nQ)\eu$ from a logarithmic derivation makes it annihilate $Q$. The arrangement is free with exponents $(1,d_1,d_2)$, ordered so that $d_1\le d_2$, precisely when $D_0(\cA)=R(-d_1)\oplus R(-d_2)$, where $d_1+d_2=n-1$ \cite{Sai80,OT92}. We also write $(d_1,d_2)$ for the two non-Euler exponents. Set
\[
r=\mdr(\cA)=\min\{k:D_0(\cA)_k\ne0\},\qquad \cE_\cA=\widetilde{D_0(\cA)}.
\]
The rank-two vector bundle $\cE_\cA$ satisfies $H^0(\cE_\cA(k))=D_0(\cA)_k$.

Throughout, $m_P$ is the number of arrangement lines through an intersection point $P$, and $s_H$ counts all intersection points on $H$, including double points. We use the convention
\begin{equation}\label{eq:comb}
b_2=\sum_P(m_P-1),\qquad
\chib(t)=t^2-(n-1)t+b_2-(n-1).
\end{equation}
Thus $b_2$ is not the constant term of the reduced characteristic polynomial $\chib$. Terao's factorisation theorem gives $\chib(t)=(t-d_1)(t-d_2)$ for a free arrangement \cite{Ter81}. A lattice with this reduced characteristic polynomial will be called of type $(n;d_1,d_2)$, whether or not it has a free realisation.

The Jacobian exact sequence
\[
0\longrightarrow\cE_\cA\longrightarrow\cO^3
\xrightarrow{(Q_x,Q_y,Q_z)}\cI_\Sigma(n-1)\longrightarrow0
\]
defines the singular scheme $\Sigma$ of $Q=0$, with $\operatorname{length}\Sigma=\sum_P(m_P-1)^2$. Indeed, at a local affine chart around a point $P$ of multiplicity $m$, the reduced local equation becomes a product $q = \ell_1 \cdots \ell_m$ of distinct linear forms, and the local Tjurina algebra is
$$
\mathcal{O}_{\mathbb{P}^2, P}/(q, q_x, q_y) \simeq \mathcal{O}_{\mathbb{P}^2, P}/(q_x, q_y),
$$
using the Euler relation. The two derivatives $q_x, q_y$ have degree $m-1$ and no common factors since $q$ is reduced with isolated singularity. Thus, they form a complete intersection at the point $P$, with multiplicity $(m-1)^2$.

Counting pairs of lines gives $\sum_Pm_P(m_P-1)=n(n-1)$, so
\begin{equation}\label{eq:chern}
c_1(\cE_\cA)=-(n-1),\qquad c_2(\cE_\cA)=b_2-(n-1)=\chib(0).
\end{equation}
We will repeatedly use the following formulas for a rank-two bundle $\cE$ on $\PP^2$:
\begin{align}
c_1(\cE(k))&=c_1(\cE)+2k,&
c_2(\cE(k))&=c_2(\cE)+kc_1(\cE)+k^2,\label{eq:twist}\\
\chi(\cE)&=2+\tfrac12c_1(\cE)(c_1(\cE)+3)-c_2(\cE),&
h^2(\cE_\cA(k))&=h^0(\cE_\cA(n-4-k)).\label{eq:RR}
\end{align}
The first equation of the second line follows using Hirzebruch--Riemann--Roch formula, while the second follows by Serre duality, using $\cE_\cA^\vee \simeq \cE_\cA \otimes \det(\cE_\cA)^{-1} = \cE_\cA(n-1)$.

A nonzero minimal syzygy $\theta$ defines a section of $\cE_\cA(r)$ with no divisorial zero: a divisor could be divided out to give a lower-degree syzygy.  The corresponding Bourbaki sequence can be written as
\begin{equation}
0 \longrightarrow \mathcal{O}(-r) \xlongrightarrow{\theta} \cE_\cA \longrightarrow \mathcal{I}_Z(r-n+1) \longrightarrow 0,
\end{equation}
defining a finite zero scheme $Z \subseteq \mathbb{P}^2$, with length
\begin{equation}\label{eq:bour}
\Bour(\cA)=c_2(\cE_\cA(r))=\chib(r).
\end{equation}
This is the Bourbaki degree \cite{NNN26}.  \cite[Theorem 2.1]{jardim2024bourbaki} establishes that length zero is equivalent to freeness. Moreover, when $r\le n/2$, length one characterises a nearly free arrangement, whose syzygy module has generators in degrees $r,n-r,n-r$, by \cite[Proposition 2.18]{jardim2024bourbaki} (see also \cite{DS18,Dim17}). We call $(r,n-r)$ its nearly free exponents. For a lattice of type $(n;d_1,d_2)$, equation~\eqref{eq:bour} reads
\begin{equation}\label{eq:bourA2}
\Bour(\cA)=(d_1-r)(d_2-r).
\end{equation}

\subsection{Restriction to a line}
For $H\in\cA$,  the distinct arrangement intersection points on $H \simeq \mathbb{P}^1$ form a reduced set $\cA^H$.  The intersections on $H$ carry multiplicities $\mu_P=m_P-1$, namely the number of lines passing through $P$ in $\cA \setminus \{H\}$, with $\sum\mu_P=n-1$. 

Choose homogeneous coordinates $[u:v]$ on $H$. A point $P$ is cut out by a non-zero linear form $\beta_P(u,v)$, unique up to multiplication by a scalar. The central rank--two multiarrangement is the collection of lines $\{\beta_P = 0: P \in \cA^H\}$ in $\mathbb{C}^2$ with multiplicities $\mu_P$, and the associated logarithmic derivation module is 
\[
D(\cA^H,\mu)=\{\eta\in\Der_\C\C[u,v]:\beta_P^{\mu_P}\mid\eta(\beta_P)\text{ for every }P\in H\}.
\]
It is free, with exponents $e_1(H)\le e_2(H)$ and $e_1+e_2=n-1$, by \cite[Corollary 7]{Zie89}. Ziegler's restriction theorem \cite[Theorem 11]{Zie89} and Yoshinaga's criterion \cite[Corollary 3.3]{Yos05} state that
\begin{equation}\label{eq:yos}
\chib(0)-e_1(H)e_2(H)\ge0,
\quad\text{with equality exactly when }\cA\text{ is free.}
\end{equation}
Consequently, for a lattice of type $(n;d_1,d_2)$,
\begin{lemma}\label{lem:yoshinaga}
For every $H$, $e_1(H)\le d_1$. The arrangement is free if and only if $e_1(H)=d_1$ for one line, equivalently for every line. This is equivalent to the absence of a nonzero restricted derivation of degree $d_1-1$.
\end{lemma}
\begin{proof}
The left side of \eqref{eq:yos} is $(e_1-d_1)(e_1-d_2)$. Since $e_1\le(d_1+d_2)/2$, nonnegativity forces $e_1\le d_1$, and equality in \eqref{eq:yos} is precisely $e_1=d_1$.
\end{proof}

If the restriction is balanced, meaning $\max\mu_P\le(n-1)/2$, Abe's bound \cite[Thm.~0.1]{Abe12} gives
\begin{equation}\label{eq:gap}
e_2(H)-e_1(H)\le s_H-2.
\end{equation}

\subsection{Candidate counterexample data}\label{sec:window}
By a counterexample we mean a free arrangement $\cA_1$ with exponents $(1,d_1,d_2)$ and a non-free arrangement $\cA_2$ with the same lattice. Write $\delta=d_2-d_1$, $r=\mdr(\cA_2)$, and $m=\max m_P$.
\begin{proposition}\label{thm:window}
A counterexample satisfies
\begin{equation}\label{eq:oldwindow}
d_1\ge6,\quad \delta\le d_1-4,\quad
\delta+4\le s_H\le d_1,\quad m\le d_1-1,\quad r\le d_1-1.
\end{equation}
The last inequality holds for every non-free arrangement of type $(n;d_1,d_2)$, without assuming a free realisation of its lattice.
\end{proposition}
\begin{proof}
The exponent bounds follow by \cite[Cor.~5.4--5.5]{Abe14}. In a free arrangement, $s_H\le d_1+1$ or $s_H=d_2+1$ \cite[Cor.~1.2]{Abe14}; equality $s_H-1\in\{d_1,d_2\}$ makes the lattice divisionally free \cite[Theorem 1.1]{Abe16}, so it cannot occur in a counterexample. The multiplicity bound follows from \cite[Cor.~1.4]{Dim17Pencils}: if the free member has $m\ge d_1$, every arrangement with its lattice is free. In particular every restriction is balanced. For $\cA_2$, Lemma~\ref{lem:yoshinaga} gives $e_1\le d_1-1$, hence $e_2-e_1\ge\delta+2$; equation~\eqref{eq:gap} gives the lower bound on $s_H$.

For the final assertion, \eqref{eq:bourA2} is positive, so $r<d_1$ or $r>d_2$. If $r>d_2$, the bundle $\cE_\cA(d_2)$ has $c_1=\delta$, $c_2=0$. Duality gives
$$
h^2(\cE_\cA(d_2))=h^0(\cE_\cA(d_1-3))=0
$$
under $r > d_2$, while Hirzebruch--Riemann--Roch gives $$
h^0(\cE_\cA(d_2)) \geq \chi(\cE_\cA(d_2)) = 2+\delta(\delta+3)/2>0,
$$ contradicting the minimality of $r>d_2$.
\end{proof}

We also need one constraint involving points on different lines. If $\cA$ is free with exponents $(1,d_1,d_2)$, Abe's external-line theorem says that a line $\ell\notin\cA$ meets it in $s=d_1+1$ or $s\ge d_2+1$ distinct points \cite[Cor.~1.2(2)]{Abe14}. Since $\sum_{P\in\ell}(m_P-1)=n-s$ (with $m_P=1$ at an ordinary crossing of $\ell$ with one arrangement line), it follows that
\begin{equation}\label{eq:unjoined}
(m_P-1)+(m_Q-1)\le d_2
\end{equation}
whenever the line $PQ$ is not an arrangement line. This is the only external-line bound needed below.

\section{A foliation bound for the least syzygy degree}\label{sec:foliation}

A logarithmic derivation is tangent to every arrangement line. After passing to the projective plane, it therefore specifies a family of curves that includes those lines. The inflection-divisor construction in Proposition~\ref{prop:infl} first gives \(n\le3r\), since its degree is \(3r\) and it contains all \(n\) arrangement lines. A comparison with the singularity count strengthens this to \(n+2\le3r\) in the non-free split-characteristic-polynomial case. This improvement is what eliminates the remaining exponent pair for seventeen lines.

\begin{theorem}\label{thm:foliation}
Let $\cA$ be a non-free arrangement of $n$ lines with
$\chib(t)=(t-d_1)(t-d_2)$, where $d_1\le d_2$. If every line contains at most $d_2+1$ intersection points, then
\[
 3\,\mdr(\cA)\ge n+2.
\]
\end{theorem}

The remainder of this section is devoted to the proof of Theorem~\ref{thm:foliation}. We first show that a minimal logarithmic derivation defines a saturated foliation of degree $r$ with isolated singularities. We then compare two global counts associated with this foliation: the degree of its inflection divisor and the total Milnor number of its singularities.

\subsection{The foliation associated with a minimal syzygy}
Take a non-zero syzygy $\theta=(a,b,c)$ of least degree $r=\mdr(\cA)$, and for the moment assume $r\ge1$ and
\begin{equation}\label{eq:hypF}
 s_H\le n-r-1\qquad(H\in\cA).
\end{equation}

The condition above is introduced here to isolate the numerical hypothesis needed below. It is automatically satisfied in the setting of Theorem~\ref{thm:foliation}. Indeed, that theorem assumes \(s_H\le d_2+1\) for every \(H\in\cA\), while Proposition~\ref{thm:window} gives \(r\le d_1-1\) for a non-free arrangement with split characteristic polynomial. Since \(d_1+d_2=n-1\), we have
\[
d_2+1=n-d_1\le n-r-1,
\]
and hence \(s_H\le n-r-1\).

Let $\Omega = dx \wedge dy \wedge dz$ be the standard volume form on $\mathbb{C}^3$ and $\eu = x \partial_x + y \partial_y + z \partial_z$ be the Euler derivation. The non-zero syzygy may be written as a vector field $\theta = a \partial_x + b \partial_y + c \partial_z$. Recall that a homogeneous vector field on the affine cone determines a projective direction only modulo the radial direction. More precisely, at a point
$p\in\mathbb C^3\setminus\{0\}$, the tangent space to projective space is
\[
T_{[p]}\PP^2\simeq T_p\mathbb C^3/\mathbb C\eu(p),
\]
so the projective line field induced by \(\theta=a\partial_x+b\partial_y+c\partial_z\) is represented upstairs by the two-plane
\[
\langle \eu(p),\theta(p)\rangle\subset T_p\mathbb C^3.
\]
This is the standard homogeneous description of a foliation on
projective space; see, for instance,
\cite[Chapter~2, \S2.2]{LinsNetoScardua2020}.

The volume form identifies such a two-plane with a one-form annihilating
it
\[
\omega=\iota_{\eu}\iota_{\theta}\Omega.
\]
Moreover, the relation
$\iota_{\eu}\omega=0$, together with homogeneity, is the Euler condition which allows $\omega$ to descend to a twisted one-form on projective
space. In coordinates one has
\[
\omega
=(bz-cy)\,dx+(cx-az)\,dy+(ay-bx)\,dz.
\]

Away from its zeros it determines a one-dimensional tangent direction on $\PP^2$. Its integral curves are the leaves of a foliation $\cF$. Because the coefficients have degree $r+1$, this foliation has degree $r$ provided that they have no common factor.

\begin{lemma}\label{lem:inv}
Let \(H\in\cA\), with \(H=\{\alpha_H=0\}\subset\PP^2.\) Along \(H\), the ambient one-form \(\omega=\iota_{\eu}\iota_{\theta}\Omega\) can be written as
\[
\omega|_H=f_H\,d\alpha_H,
\]
where \(f_H\) is a non-zero binary form of degree \(r+1\).
Its zeros are precisely the singularities of \(\cF\) lying on \(H\),
and every intersection point of the arrangement lying on \(H\) is among
these zeros. Writing \(\nu_P(H)=\ord_P(f_H),\)
one has
\begin{equation}\label{eq:line-zero-count}
\sum_{P\in H}\nu_P(H)=r+1.
\end{equation}
\end{lemma}

\begin{proof}
 Since
\(\theta\in D_0(\cA)\subset D(\cA)\),  \(\alpha_H\mid\theta(\alpha_H),\) so there is a homogeneous polynomial \(K_H\) of degree \(r-1\) such that \(\theta(\alpha_H)=K_H\alpha_H.\)
Since \(\eu(\alpha_H)=\alpha_H\), 
we have \(\eu(\alpha_H)=0\) on \(H\).

After a linear change of coordinates we may assume that
\(\alpha_H=z\), so \(H=\{z=0\}\). The logarithmicity condition along \(H\) gives
\[
c=\theta(z)=K_Hz,
\]
hence \(c|_H=0\). Since
\[
\omega=\iota_{\eu}\iota_{\theta}
       (dx\wedge dy\wedge dz)
=(bz-cy)\,dx+(cx-az)\,dy+(ay-bx)\,dz,
\]
restriction to \(H\) gives \(\omega|_H=(ay-bx)|_{z=0}\,dz.\)
Thus \(\omega|_H=f_H\,d\alpha_H\) for a homogeneous binary form $f_H$ of degree $r+1$, which measures the linear dependency of the two tangent vectors $\theta|_H$ and $\eu|_H$ on the affine cone over $H$.

Moreover, let
\(P=H\cap H'\) be an intersection with another arrangement line
\(H'\neq H\). At a representative \(p\in\mathbb C^3\setminus\{0\}\) of
\(P\), logarithmicity implies that \(\theta(p)\) is tangent to both
affine cones over \(H\) and \(H'\). The intersection of these two
tangent planes is the radial line \(\mathbb Cp\), so \(\theta(p)\in\mathbb C\,\eu(p).\) Hence the two vectors contracted into \(\Omega\) are linearly dependent
at \(p\), and therefore
\[
\omega(p)=\iota_{\eu(p)}\iota_{\theta(p)}\Omega=0.
\]
Thus every intersection point of the arrangement on \(H\) is a zero of
\(f_H\).

It remains to prove that \(f_H\) is not identically zero. Let \(\theta_H=\theta-K_H\eu\), which is also logarithmic along $\cA$, and
\[
\theta_H(\alpha_H)
=\theta(\alpha_H)-K_H\eu(\alpha_H)
=K_H\alpha_H-K_H\alpha_H
=0.
\]
Hence \(\theta_H\) belongs to the submodule of logarithmic derivations
annihilating \(\alpha_H\). By the standard Ziegler multirestriction construction, the restriction $\theta_H|_H$ belongs to \(D(\cA^H,\mu)\); see the proof of \cite[Example~2 and the proof of Theorem~11]{Zie89}. We note that replacing \(\theta\) by \(\theta_H\) does not change the one-form
defining the foliation, since
\[
\iota_{\eu}\iota_{\theta_H}\Omega
=
\iota_{\eu}\iota_{\theta-K_H\eu}\Omega
=
\iota_{\eu}\iota_{\theta}\Omega
-
K_H\,\iota_{\eu}\iota_{\eu}\Omega
=
\omega,
\]
using \(\iota_{\eu}^2=0\). In particular, the same binary form \(f_H\)
is obtained from \(\theta_H\).

Suppose now that \(f_H\equiv0\). In the coordinates \(H=\{z=0\}\), write
\[
\theta_H|_H=A(x,y)\partial_x+B(x,y)\partial_y.
\]
Since \(\theta_H(\alpha_H)=0\), there is no normal component, and from
the computation above we have
\[
f_H=A(x,y)y-B(x,y)x.
\]
Thus \(f_H\equiv0\) implies \(Ay-Bx=0.\) Since \(x\) and \(y\) are relatively prime in \(\mathbb C[x,y]\), there
exists a homogeneous polynomial \(g\) of degree \(r-1\) such that
\[
A=gx,\qquad B=gy.
\]
Equivalently, \(\theta_H|_H=g\,\eu_H,\) where \(\eu_H=x\partial_x+y\partial_y\) is the Euler derivation on the affine cone over \(H\simeq\PP^1\).

If \(g=0\), then every coefficient of \(\theta_H\) is divisible by
\(\alpha_H\), so
\[
\theta_H=\alpha_H\eta
\]
for a homogeneous derivation \(\eta\) of degree \(r-1\). Since
\(\theta_H\) is logarithmic along every line of \(\cA\), and
\(\alpha_H\) is relatively prime to the defining equation of any other
line, it follows that \(\eta\in D(\cA)\). Moreover, \(\eta\) cannot be
an Euler multiple: otherwise \(\theta=K_H\eu+\alpha_H\eta\)
would itself be an Euler multiple, contradicting \(0\neq\theta\in D_0(\cA)\). Therefore the non-Euler component of
\(\eta\) gives a non-zero element of \(D_0(\cA)\) of degree \(r-1\),
contradicting the minimality of
\(r=\mdr(\cA)\).

Assume instead that \(g\neq0\). Since
\(\theta_H|_H\in D(\cA^H,\mu)\), for every intersection point
\(P\in H\) we have
\[
\beta_P^{\mu_P}
 \mid (\theta_H|_H(\beta_P))
 =g\,\eu_H(\beta_P)
 =g\,\beta_P.
\]
Consequently, \(\beta_P^{\mu_P-1}\mid g\) for every $P\in \cA^H.$
The forms \(\beta_P\) corresponding to distinct points of \(H\) are
pairwise relatively prime, and therefore
\[
r-1=\deg g
\ge\sum_{P\in H}(\mu_P-1)
=\sum_{P\in H}\mu_P-s_H
=n-1-s_H.
\]
By \eqref{eq:hypF}, \(s_H\le n-r-1,\) so \(n-1-s_H\ge r\). Hence \(r-1\ge r\), a contradiction.

Thus \(f_H\not\equiv0\). Since \(f_H\) is a non-zero binary form of
degree \(r+1\) on \(H\simeq\PP^1\), the sum of the orders of its zeros
is \(r+1\), proving \eqref{eq:line-zero-count}.
\end{proof}

\begin{lemma}\label{lem:sat}
The coefficients of $\omega$ have no common factor. Therefore, $\cF$ is a foliation of degree $r$ with isolated singularities. Moreover, their Milnor numbers satisfy
\begin{equation}\label{eq:total-milnor}
\sum_{P\in\Sing\cF}\mu_P=r^2+r+1.
\end{equation}
\end{lemma}

\begin{proof}
 Suppose, by contradiction, that \(\omega=h \omega'\) for a homogeneous polynomial \(h\) of degree \(k>0\) and a homogeneous one-form \(\omega'\). Since \(0=\iota_{\eu}\omega=h\,\iota_{\eu}\omega'\) and \(h\neq0\), we also have \(\iota_{\eu}\omega'=0.\) Writing
\[
\omega'=A\,dx+B\,dy+C\,dz,
\]
this condition is \(xA+yB+zC=0.\) By exactness of the Koszul complex associated with \(x,y,z\), there are homogeneous polynomials $a', b', c'$ such that
$$
(A,B,C) = (b'z - c' y, c'x -a'z, a' y - b'x).
$$
Setting $\theta' = a' \partial_x + b' \partial_y + c' \partial_z$, the identity above reads  \(\omega'=\iota_{\eu}\iota_{\theta'}\Omega.\)
Therefore
\[
\omega
=h\,\iota_{\eu}\iota_{\theta'}\Omega
=\iota_{\eu}\iota_{h\theta'}\Omega \Longrightarrow \iota_{\eu}\iota_{\theta-h\theta'}\Omega=0.
\]
By the Euler sequence, a homogeneous vector field on the affine cone induces the zero vector field on $\mathbb{P}^2$ precisely when it is a multiple of the Euler vector field. Hence \(\theta=h\theta'+g\eu\) for some homogeneous polynomial \(g\).

We claim that \(\theta'\) is logarithmic along every line of \(\cA\).
Fix \(H=\{\alpha_H=0\}\in\cA\). By Lemma~\ref{lem:inv},  \(\omega|_H=f_H\,d\alpha_H\),
with \(f_H\not\equiv0\). In particular, \(\omega\) does not vanish identically along \(H\). If \(\alpha_H\mid h\), then \(h|_H=0\), and the factorization \(\omega=h\omega'\) would imply \(\omega|_H=0\), a contradiction.

On the other hand, since \(\theta\in D(\cA)\), \(\alpha_H\mid\theta(\alpha_H).\) Using \(\theta=h\theta'+g\eu\)
and \(\eu(\alpha_H)=\alpha_H\), we obtain
\[
\theta(\alpha_H)
=
h\,\theta'(\alpha_H)+g\alpha_H.
\]
Reducing modulo \(\alpha_H\) gives \(h\,\theta'(\alpha_H)\equiv0\pmod{\alpha_H}.\) Since \(\alpha_H\) is irreducible and \(\alpha_H\nmid h\), it follows
that \(\alpha_H\mid\theta'(\alpha_H)\), so \(\theta'\) is logarithmic along each \(H\), giving $\theta' \in D(\cA)$.

Moreover, \(\theta'\) is not an Euler multiple. Indeed, if
\(\theta'=q\eu\), then \(\omega'
=
\iota_{\eu}\iota_{q\eu}\Omega
=0,\) contrary to \(\omega=h\omega'\neq0\). Hence, using the decomposition \(D(\cA)=R\eu\oplus D_0(\cA),\)
the \(D_0(\cA)\)-component of \(\theta'\) is non-zero. Its degree is
\[
\deg\theta'=r-k<r,
\]
contradicting the minimality of
\(r=\mdr(\cA)\). Therefore the coefficients of \(\omega\) have no common factor.

The homogeneous one-form \(\omega\) defines a section \(\omega\in H^0\bigl(\PP^2,\Omega^1_{\PP^2}(r+2)\bigr)\). Since its coefficients have no common factor, its zero scheme has no divisorial component and is therefore zero-dimensional. These zeros are precisely the singularities of the foliation \(\cF\).

To identify their local lengths, let $P$ be a singular point. We may assume that $P \in U_z = \{z \neq 0\} \subseteq \mathbb{P}^2$, after changing coordinates. Thus, we may work on the affine chart $z = 1$, with coordinates $u = x/z$ and $v = y/z$. Writing \(\theta=a\partial_x+b\partial_y+c\partial_z\), the restriction of \(\omega\) to this chart is
\[
\omega_{\mathrm{aff}}
=
(b-vc)\,du+(uc-a)\,dv.
\]
Setting \(X_1=a-uc\) and \(X_2=b-vc\), where \(a,b,c\) are evaluated at \((u,v,1)\), we obtain
\begin{equation}\label{eq:foliation-affine-expression}
\omega_{\mathrm{aff}}
=
X_2\,du-X_1\,dv,    
\end{equation}
so the local zero scheme of \(\omega\) at a singular point \(P\) is defined by the ideal \((X_1,X_2)\subset\cO_P.\) Accordingly, its local length is
\[
\dim_{\C}\frac{\cO_P}{(X_1,X_2)}
=\mu_P,
\]
the Milnor number of the isolated singularity of \(\cF\) at \(P\) (see, for example, \cite[Definition~6]{Alcantara2011}). The global Milnor-number formula for a degree-\(r\) foliation on \(\PP^2\) with isolated singularities now gives
\[
\sum_{P\in\Sing\cF}\mu_P
=
r^2+r+1;
\]
see \cite[Proposition~2]{Alcantara2011}. This proves
\eqref{eq:total-milnor}.
\end{proof}

\subsection{The inflection divisor and its local multiplicities}

The inflection condition is local. Let $P$ be a regular point of $\cF$ and choose an affine chart containing $P$, say $z \neq 0$. Using the same notation $u = x/z$, $v = y/z$ as in the proof of Lemma~\ref{lem:sat}, we may consider the vector field $X = X_1 \partial_u + X_2 \partial_v$ such that
$$
\iota_X(du \wedge dv) = X_1 dv - X_2 du = - \omega_{\mathrm{aff}}
$$
as in \eqref{eq:foliation-affine-expression}. The vector field $X$ defines the same foliation as $\omega$ in this affine chart.

At a local analytic neighbourhood of $P$, choose a local holomorphic integral curve $\gamma$, parametrized by its flow parameter, so that $\gamma' = X(\gamma)$. Hence $\gamma'' = DX(\gamma) X(\gamma)$, where
$$
DX = \begin{pmatrix}
\partial_u X_1 & \partial_v X_1\\
\partial_u X_2 & \partial_v X_2
\end{pmatrix}
$$
is the associated Jacobian matrix.

The tangent line of $\cF$ at $P$ has contact of order at least three if and only if $\det(X, DX \cdot X) = 0$. The corresponding global equation on \(\PP^2\) is the \emph{first extactic determinant}; see \cite[\S3.2]{Per01}. In terms of the homogeneous lift $\theta = a \partial_x + b \partial_y + c \partial_z$, it is the equation 
\begin{equation}\label{eq:extactic}
 J=\det\begin{pmatrix}
 x&y&z\\ a&b&c\\ \theta(a)&\theta(b)&\theta(c)
 \end{pmatrix}.
\end{equation}
The three rows have degrees $1$, $r$ and $2r-1$, respectively, and $J$ defines the \emph{first extactic curve} $\mathcal{E}_1(\omega)$.

\begin{proposition}\label{prop:infl}
The polynomial $J$ is non-zero, has degree $3r$, and is divisible by $Q$, so in particular $3r \geq n$. Moreover, $J$ vanishes to order at least two at every singular point of $\cF$.
\end{proposition}

\begin{proof}
Suppose that $J\equiv0$. \cite[Thm.~2]{Per01} shows that the vanishing of the $(d)$-th extactic curve $\mathcal{E}_d(\omega)$ forces the foliation to admit a meromorphic first integral of degree at most $(d)$. Hence \(J\equiv0\) implies that \(\cF\) admits a first integral of degree one. Its fibers are therefore the lines of a pencil, so the saturated foliation has degree zero. This contradicts \(r\ge1\) and Lemma~\ref{lem:sat}. Thus \(J\not\equiv0\). The corresponding degree assertion follows from the degrees of the rows in \eqref{eq:extactic}. 

We next prove that the arrangement is contained in the extactic divisor.
For every \(H=\{\alpha_H=0\}\in\cA\), logarithmicity gives \(\theta(\alpha_H)=K_H\alpha_H,\) so \(H\) is an invariant line of the foliation \(\cF\). Since \(J\) is
the first extactic polynomial \(\mathcal E_1(\cF)\), \cite[Proposition~1]{Per01} implies that the equation of every
invariant line is a factor of \(J\). Therefore
\[
\alpha_H\mid J
\qquad\text{for every }H\in\cA.
\]
The lines of \(\cA\) are distinct, so their defining linear forms are
pairwise non-associated irreducible factors. Consequently,
\[
Q=\prod_{H\in\cA}\alpha_H\mid J.
\]

Finally, let \(P\in\Sing\cF\). Choose an affine chart containing \(P\)
and local coordinates \((u,v)\) centered at \(P\), as above, with corresponding vector field \(X=X_1\partial_u+X_2\partial_v\). Then
\[
J_{\mathrm{aff}}
=
\det\bigl(X,DX\cdot X\bigr).
\]
Since \(P\) is a singular point, \(X_1(P)=X_2(P)=0\), and hence \(X_1,X_2\in\mathfrak m_P,\) where \(\mathfrak m_P\) is the maximal ideal of the local ring at \(P\).
Moreover, each component of \(DX\cdot X\) has the form
\[
(\partial_u X_i) X_1
+ (\partial_v X_i) X_2,
\]
and therefore also belongs to \(\mathfrak m_P\). It follows that \(J_{\mathrm{aff}} \in\mathfrak m_P^2\), and thus \(\mult_P(J)\ge2\)
for every \(P\in\Sing\cF\), as claimed.
\end{proof}

Write $J=QJ'$ and let $j = \deg J'=3r-n\ge0$. To prove the desired bound \(3r\ge n+2\) of Theorem~\ref{thm:foliation}, it remains to rule out the two possibilities $\deg(J') = 0$ or $\deg(J') = 1$ when $\cA$ is not free. The next lemma provides the local estimate needed: it compares the contribution of an arrangement singularity to the Milnor number of the foliation at each point. Summing these local contributions will then be combined with the global Milnor-number formula to exclude these cases.

\begin{lemma}\label{lem:degenerate}
At a singular point $P$, let $L_P$ be the linear part of $X$, and let $\ell_P$ be the number of invariant projective lines through $P$. If $\det L_P=0$, then $\mult_P(J)\ge3$. If $L_P=0$, then
\[
 \mult_P(J)\ge\ell_P+2.
\]
\end{lemma}

\begin{proof}
Consider again the affine coordinates $(u,v)$ centered at $P$. Explicitly, one has
$$
X(u,v) = \begin{pmatrix}
X_1(u,v)\\
X_2(u,v)
\end{pmatrix} = L_P \begin{pmatrix}
u\\
v
\end{pmatrix} + \text{ terms of degree }\geq 2, \qquad L_P = \begin{pmatrix}
\partial_u X_1(P) & \partial_v X_1(P)\\
\partial_u X_2(P) & \partial_v X_2(P).
\end{pmatrix}
$$
The quadratic term of
$\det(X,DX\cdot X)$ is $\det(L_Pt,L_P^2t)$, denoting $t = (u,v)$. By the Cayley--Hamilton identity, \(L_P^2-(\operatorname{tr}L_P)L_P+(\det L_P)\operatorname{Id}=0.\) If \(\det L_P=0\), then \(L_P^2=(\operatorname{tr}L_P)L_P\), and thus the quadratic
part of \(J_{\mathrm{aff}}\) vanishes, giving \(\mult_P(J)\ge3.\)

Suppose now that \(L_P=0\). Write the homogeneous expansion of the local
vector field as
\[
X=X_k+X_{k+1}+\cdots,
\qquad k\geq2,
\]
where each \(X_a\) is homogeneous of degree \(a\). Let \(\eu_0=u\partial_u+v\partial_v\) denote the radial vector field, and let \(i\geq k\) be the smallest
index such that \(X_i\) is not a polynomial multiple of \(\eu_0\). If every homogeneous component of \(X\)
were radial, then \(X=h\eu_0\) for some polynomial \(h\), so
\(DX\cdot X\) would also be radial. It would follow that
\[
J_{\mathrm{aff}}
=
\det(X,DX\cdot X) = \det( h\eu_0, (h + \eu_0(h))h\eu_0 )
\equiv0,
\]
contrary to the fact that \(J\not\equiv0\). 

Let \(L\) be an invariant projective line through \(P\). In the chosen
affine coordinates, \(L\) is a line through the origin. Since \(X\) is
tangent to \(L\), the vectors \(t\) and \(X(t)\) are linearly dependent
at every point of \(L\), and hence \(\det(t,X(t))=0\)
on $L$, so the equation of $L$ divides $\det(t,X(t))$. For every \(a<i\), the field \(X_a\) is radial, hence $X_a = q_a \eu_0$, so at $t= (u,v)$:
\[
\det(t,X_a) =\det(t, q_a(t)t) = q_a(t) \det(t,t) =0.
\]
Therefore the first non-zero homogeneous term of \(\det(t,X)\) is \(\det(t,X_i).\)
Since \(X_i\) is not radial, this is a non-zero binary form of degree
\(i+1\). Each invariant line through \(P\) gives a distinct linear factor
of this binary form. Consequently, \(\ell_P\leq i+1.\)

We now expand
\[
J_{\mathrm{aff}}
=
\sum_{a,b,c}
\det\bigl(X_c,DX_a\cdot X_b\bigr),
\]
so the summand indexed by \((a,b,c)\) is homogeneous of degree \(a+b+c-1.\) If \(a,b,c<i\), then \(X_a,X_b,X_c\) are all radial, and $  \det\bigl(X_c,DX_a\cdot X_b\bigr)=0.$ Thus every non-zero term in the expansion of \(J_{\mathrm{aff}}\) has
at least one of \(a,b,c\) greater than or equal to \(i\). Since all three indices are at least \(k\),
\[
a+b+c-1\geq 2k+i-1.
\]
As \(k\geq2\), \(2k+i-1\geq i+3,\) and, since \(\ell_P\leq i+1\), \(i+3\geq\ell_P+2\). Therefore every non-zero homogeneous term of \(J_{\mathrm{aff}}\) has
degree at least \(\ell_P+2\). Hence \(\mult_P(J)\geq\ell_P+2\), showing the claim. 
\end{proof}

\subsection{The excess identity and proof of Theorem~\ref{thm:foliation}}

For a singular point $P$ of $\cF$, let $m'_P$ be the number of arrangement lines through it, including the possibility $m'_P=0$. For points on the arrangement define
\[
 d_P=\sum_{H\ni P}\nu_P(H)-(m'_P-1)-\mu_P.
\]
The quantity $d_P$ measures the discrepancy between counting this singularity on its invariant lines and counting its local scheme length.

\begin{lemma}\label{lem:local}\label{lem:tame}
Let $P$ be a singular point of $\cF$ on the arrangement, and let $L_P = DX(P)$ be the linearization of the local vector field $X$ defining $\cF$ at $P$. Then $d_P = 0$ if either one of the following cases hold.
\begin{itemize}
    \item[(a)] $\det(L_P) \neq 0$;
    \item[(b)] \(m'_P\in\{1,2\}\) and, in coordinates adapted to an arrangement line \(H\) through \(P\),
    \(L_P=\operatorname{diag}(0,e)\) with \(e\neq0,\)
    and \(H\) corresponding to the zero eigendirection.
\end{itemize}
\end{lemma}

\begin{proof}
Choose affine coordinates \((u,v)\) centred at \(P\), and write \(X=X_1(u,v)\partial_u+X_2(u,v)\partial_v\). If \(H\) is an invariant line through \(P\), with direction generated by a non-zero vector \(\xi\), then \(X(t\xi)\in\C\xi\) for all sufficiently small \(t\). Since
\[
X(t\xi)=tL_P\xi+O(t^2),
\]
it follows that \(L_P\xi\in\C\xi\). Thus every invariant line through
\(P\) determines an eigendirection of \(L_P\).

Suppose that \(H\) is given by \(v=0\). Since \(H\) is invariant, \(X_2(u,0)=0\), and the induced one-dimensional vector field on \(H\) is
\[
X|_H=X_1(u,0)\partial_u.
\]
Moreover, from \(\omega_{\mathrm{aff}}=X_2\,du-X_1\,dv\) we obtain \(\omega|_H=-X_1(u,0)\,dv.\) Hence
\[
\nu_P(H)=\ord_0X_1(u,0).
\]
In particular, if the eigenvalue of \(L_P\) corresponding to \(H\) is
non-zero, then the linear term of \(X_1(u,0)\) is non-zero, and therefore \(\nu_P(H)=1.\)

Assume first that \(\det L_P\neq0\). Then every eigenvalue of \(L_P\) is
non-zero, so \(\nu_P(H)=1\) for every arrangement line \(H\) through \(P\). Consequently,
\[
\sum_{H\ni P}\nu_P(H)=m'_P.
\]
Since \(L_P = DX(P)\) is invertible, the map $(u,v) \mapsto (X_1, X_2)$ has invertible Jacobian at $P$. By the holomorphic inverse function theorem, $(X_1, X_2)$ form a system of local analytic coordinates at $P$, so $(X_1, X_2) = \mathfrak{m}_P$ in the local ring. Therefore
\[
\mu_P
=
\dim_\C\frac{\cO_P}{(X_1,X_2)}
=1.
\]
It follows that \(d_P = m'_P-(m'_P-1)-1 =0.\)

Assume now that \(L_P=\operatorname{diag}(0,e)\), with $e\neq0$, and that \(H=\{v=0\}\) is the zero eigendirection. Since \(H\) is invariant, $X_2(u,0) = 0$, and we may then write \(X_2=vq\) for some local holomorphic function \(q\) around $P$. Since the linear part of \(X_2\) is \(ev\), we have \(q(P)=e\neq0\). Thus \(q\) is a unit in
\(\cO_P\), and hence \((X_1,X_2)=(X_1,v)\) as ideals in the local ring. Therefore
\[
\mu_P
=
\dim_\C\frac{\cO_P}{(X_1,v)}
=
\ord_0X_1(u,0)
=
\nu_P(H).
\]

If \(m'_P=1\), then \(H\) is the only arrangement line through \(P\), and hence
\[
d_P
=
\nu_P(H)-(1-1)-\mu_P
=0.
\]
If \(m'_P=2\), let \(K\) be the second arrangement line through \(P\).
Its direction is another eigendirection of \(L_P\), hence it must
correspond to the non-zero eigenvalue \(e\). Therefore \(\nu_P(K)=1\)
and
\[
d_P
=
\nu_P(H)+\nu_P(K)-(2-1)-\mu_P
=
\mu_P+1-1-\mu_P
=0.
\]
\end{proof} 

\begin{proposition}[excess identity]\label{thm:excess}
Under \eqref{eq:hypF}, if $\chib(t)=(t-d_1)(t-d_2)$, then
\begin{equation}\label{eq:excess}
 (d_1-r)(d_2-r)
 =\sum_{P\notin\bigcup\cA}\mu_P
  -\sum_{P\in\bigcup\cA}d_P,
\end{equation}
where both sums are over singularities of $\cF$.
\end{proposition}

\begin{proof}
By Lemma~\ref{lem:inv}, every arrangement line \(H\in\cA\) satisfies
\[
\sum_{P\in H}\nu_P(H)=r+1.
\]
Summing over the \(n\) lines of \(\cA\) and then regrouping according to
the singular points of \(\cF\) lying on the arrangement gives
\[
n(r+1)
=
\sum_{P\in\Sing\cF\cap\bigcup\cA} \left(
\sum_{\substack{H\in\cA\\P\in H}}\nu_P(H) \right).
\]
By the definition of \(d_P\),
\[
\sum_{H\ni P}\nu_P(H)
=
(m'_P-1)+\mu_P+d_P,
\]
and therefore
\[
n(r+1)
=
\sum_{P\in\Sing\cF\cap\bigcup\cA}(m'_P-1)
+\sum_{P\in\Sing\cF\cap\bigcup\cA}\mu_P+\sum_{P\in\Sing\cF\cap\bigcup\cA}d_P.
\]

The first sum is \(b_2\). Indeed, every intersection point of the
arrangement is a singular point of \(\cF\) by Lemma~\ref{lem:inv},
and at such a point \(m'_P=m_P\); a singularity of \(\cF\) lying on
only one arrangement line has \(m'_P=1\) and hence contributes zero.
Thus, by \eqref{eq:comb},
\[
\sum_{P\in\Sing\cF\cap\bigcup\cA}(m'_P-1)
=
\sum_{P\in\Sing\cA}(m_P-1)
=
b_2.
\]

On the other hand, Lemma~\ref{lem:sat} gives
\[
\sum_{P\in\Sing\cF}\mu_P=r^2+r+1.
\]
Hence
\[
\sum_{P\in\Sing\cF\cap\bigcup\cA}\mu_P
=
r^2+r+1
-
\sum_{P\in\Sing\cF\setminus\bigcup\cA}\mu_P.
\]
Substituting these two identities yields
\[
n(r+1)
=
b_2+r^2+r+1
-\sum_{P\in\Sing\cF\setminus\bigcup\cA}\mu_P
+\sum_{P\in\Sing\cF\cap\bigcup\cA}d_P.
\]

Finally, since
\[
\chib(t)
=
t^2-(n-1)t+b_2-(n-1)
=
(t-d_1)(t-d_2),
\]
we have \(d_1+d_2=n-1\) and $b_2=d_1d_2+n-1.$ Consequently,
\[
\begin{aligned}
b_2+r^2+r+1-n(r+1)
&=d_1d_2+r^2-r(n-1)\\
&=d_1d_2+r^2-r(d_1+d_2)\\
&=(d_1-r)(d_2-r).
\end{aligned}
\]
Rearranging proves \eqref{eq:excess}.
\end{proof}

\begin{proof}[Proof of Theorem~\ref{thm:foliation}]
Set \(r=\mdr(\cA)\). By the final assertion of
Proposition~\ref{thm:window}, every non-free arrangement with \(\chib(t)=(t-d_1)(t-d_2)\) satisfies \(r\leq d_1-1.\) Since \(d_1+d_2=n-1\), the hypothesis of the theorem gives, for every \(H\in\cA\),
\[
s_H\leq d_2+1=n-d_1\leq n-r-1.
\]
Thus \eqref{eq:hypF} holds, and all the preceding results of this section apply.

We also have \(r\geq1\). Indeed, if \(r=0\), then there is a non-zero
constant Jacobian syzygy. After a linear change of coordinates, we may
assume that it is \(\partial_z\), so \(Q_z=0\). Hence \(Q\in\C[x,y]\),
and all the lines of \(\cA\) pass through a common point. Thus \(\cA\)
is a pencil and is free, contrary to hypothesis.

Moreover, \(d_1-r\geq1\), and $d_2-r\geq d_1-r\geq1$, and therefore
\begin{equation}\label{eq:positive-excess}
(d_1-r)(d_2-r)>0.
\end{equation}

Recall from Proposition~\ref{prop:infl} that
\[
J=QJ',
\qquad
j:=\deg J'=3r-n.
\]
Suppose, by contradiction, that \(j\leq1\). Since \(Q\) is the reduced
product of the arrangement lines, its multiplicity at a point \(P\) is
precisely the number \(m'_P\) of arrangement lines through \(P\). Hence
\[
\mult_P(J)
=
m'_P+\mult_P(J').
\]
Since \(J'\) has degree at most one,
\(\mult_P(J')\leq1\) at every point, and therefore
\begin{equation}\label{eq:small-residual}
\mult_P(J)\leq m'_P+1.
\end{equation}

We claim that every singularity of \(\cF\) lies on the arrangement. Assume otherwise, so \(P\in\Sing\cF\) and
\(P\notin\bigcup\cA\). Then \(m'_P=0\), so \eqref{eq:small-residual} gives \(\mult_P(J)\leq1.\) On the other hand, Proposition~\ref{prop:infl} gives
\[
\mult_P(J)\geq2
\]
at every singular point of \(\cF\), a contradiction. Consequently, \(\Sing\cF\subseteq\bigcup\cA\).

Let now \(P\in\Sing\cF\), and let \(L_P=DX(P)\) be the linearization of
a local vector field defining \(\cF\) at \(P\). We classify the possibilities of $L_P$.

First, if \(L_P=0\), then Lemma~\ref{lem:degenerate} gives \(\mult_P(J)\geq\ell_P+2,\) where \(\ell_P\) is the number of invariant projective lines through
\(P\). Since every arrangement line is invariant,
\(\ell_P\geq m'_P\), and
\[
\mult_P(J)\geq m'_P+2,
\]
contradicting \eqref{eq:small-residual}. Thus $L_P \neq 0$.

Suppose next that \(L_P\neq0\) but \(\det L_P=0\).
Lemma~\ref{lem:degenerate} then gives \(\mult_P(J)\geq3.\) Together with \eqref{eq:small-residual}, this implies \(m'_P\geq2.\)
Every arrangement line through \(P\) determines an eigendirection of \(L_P\), as in the proof of Lemma~\ref{lem:tame}. Since \(L_P\) is non-zero and singular, it is not scalar: a scalar singular matrix would necessarily be zero. A non-scalar endomorphism of a two-dimensional vector space has at most two distinct eigendirections. Hence \(m'_P=2.\) The two arrangement lines through \(P\) give two distinct
eigendirections, so \(L_P\) is diagonalizable. Since
\(\det L_P=0\) and \(L_P\neq0\), in coordinates adapted to these two lines we have
\[
L_P=\operatorname{diag}(0,e),
\qquad e\neq0.
\]
Thus case~(2) of Lemma~\ref{lem:tame} applies and gives \(d_P=0.\)

If instead \(\det L_P\neq0\), then case~(1) of
Lemma~\ref{lem:tame} gives again \(d_P=0.\)
We conclude that every singular point of \(\cF\) lies on the arrangement and satisfies \(d_P=0\).

Applying the excess identity \eqref{eq:excess}, the sum over singularities outside the arrangement is empty and the discrepancy term vanishes. Hence \((d_1-r)(d_2-r)=0,\) contradicting \eqref{eq:positive-excess}. Therefore \(j\geq2\), that is, \(3r-n\geq2\), as claimed.
\end{proof}

\begin{corollary}\label{cor:foliation-numerical-consequence}
Suppose that a counterexample has type \((n;d_1,d_2)\), with \(d_1\leq d_2\). Then
\[
d_2-d_1\leq d_1-6.
\]
In particular, \(d_1\geq7\) and \(n\geq15\).
\end{corollary}

\begin{proof}
For a counterexample, Proposition~\ref{thm:window} gives \(r\leq d_1-1\), while Theorem~\ref{thm:foliation} gives \(3r\geq n+2.\) Therefore
\[
3(d_1-1)\geq n+2.
\]
Since \(d_1+d_2=n-1\), equivalently \(n=d_1+d_2+1\), this yields
\[
3d_1-3\geq d_1+d_2+3,
\]
and hence \(d_2-d_1\leq d_1-6.\)

If \(d_1=6\), then \(d_2\geq d_1\) together with the preceding
inequality forces \(d_2=6\), so
\[
n=d_1+d_2+1=13.
\]
The thirteen-line case is already known by \cite{DIM19}, so no
counterexample can occur. Hence \(d_1\geq7\). Since \(d_2\geq d_1\),
we also obtain
\[
n=d_1+d_2+1\geq2d_1+1\geq15.
\]
\end{proof}

\section{The seventeen-line reduction}\label{sec:seventeen}

The foliation bound is useful because its integer rounding eliminates one of the two remaining exponent pairs.
\begin{theorem}\label{thm:17}
In a seventeen-line counterexample the free exponents are $(1,8,8)$, and the non-free member is nearly free with exponents $(7,10)$. Its lattice satisfies
\begin{enumerate}[label=(W\arabic*),ref=W\arabic*]
\item $b_2=80$, or equivalently $\chib(t)=(t-8)^2$;
\item $m_P\le6$ at every point;
\item $4\le s_H\le8$ on every line;
\item every eight-point line contains a double point;
\item $(m_P-1)+(m_Q-1)\le8$ whenever $P,Q$ are not joined by an arrangement line.
\end{enumerate}
More generally, every non-free seventeen-line arrangement with $\chib(t)=(t-8)^2$ and $s_H\le9$ is nearly free with exponents $(7,10)$.
\end{theorem}
\begin{proof}
Let \(\cA_1\) be the free member of the counterexample and \(\cA_2\) the non-free member, and write the exponents of \(\cA_1\) as \((1,d_1,d_2)\), with \(d_1\leq d_2\). Since \(d_1+d_2=n-1=16,\) Corollary~\ref{cor:foliation-numerical-consequence} gives \(d_2-d_1\leq d_1-6.\) Substituting \(d_2=16-d_1\), we obtain \(d_1\geq8\). Since \(d_1\leq d_2\) and
\(d_1+d_2=16\), it follows that \(d_1=d_2=8.\)
Thus the free member has exponents \((1,8,8)\), and
\[
\chib(t)=(t-8)^2.
\]

We next prove the more general assertion in the statement. Let \(\cA\)
be any non-free arrangement of seventeen lines satisfying
\[
\chib(t)=(t-8)^2
\qquad\text{and}\qquad
s_H\leq9
\quad\text{for every }H\in\cA.
\]
The final assertion of Proposition~\ref{thm:window}, which applies to
every non-free arrangement of type \((17;8,8)\), gives \(r=\mdr(\cA)\leq7.\) On the other hand, the hypothesis \(s_H\leq9=d_2+1\) allows us to
apply Theorem~\ref{thm:foliation}, yielding \(3r\geq19.\) Hence \(r\geq7\), and therefore $r = 7$.

Equation~\eqref{eq:bourA2} now gives \(\Bour(\cA)=(8-7)^2=1.\) Since \(r=7<17/2\), the length-one characterization recalled after
\eqref{eq:bour} shows that \(\cA\) is nearly free, with exponents \((r,17-r)=(7,10)\). This applies to the non-free member \(\cA_2\), since Proposition~\ref{thm:window} gives \(s_H\leq8\) for every \(H\in\cA_2\).

We now establish \textup{(W1)--(W5)}. Since
\(\chib(t)=(t-8)^2\), comparison with \eqref{eq:comb} gives \(b_2=80,\) proving \textup{(W1)}. Moreover, Proposition~\ref{thm:window}, with
\(\delta=d_2-d_1=0\), gives \(4\leq s_H\leq8\) for every line \(H\), proving \textup{(W3)}. Finally,
\eqref{eq:unjoined}, with \(d_2=8\), gives
\[
(m_P-1)+(m_Q-1)\leq8
\]
whenever \(P\) and \(Q\) are not joined by an arrangement line, which
is \textup{(W5)}.

For \textup{(W2)}, put $m=\max_P m_P$. Proposition~\ref{thm:window} gives $m\le7$. The deduction of \textup{(W2)} uses the near-freeness criterion in the form stated in \cite[Prop.~3.10]{AIM24}, following \cite[Cor.~1.7]{Dim17Pencils}: a nearly free arrangement with exponents $(a,b)$ and $m\ge a$ has only nearly free realisations of its lattice. Under this criterion, $m=7$ forces every realisation of the lattice of $\cA_2$ to be nearly free, since its exponents are $(7,10)$. This contradicts the freeness of $\cA_1$ and gives $m\le6$.

Finally, suppose that an eight-point line \(H\) contains no double
point. The other sixteen arrangement lines meet \(H\), and so
\[
\sum_{P\in H}(m_P-1)=16.
\]
Since \(H\) contains exactly eight intersection points and none is double, each term satisfies \(m_P-1\geq2\). Therefore equality forces \(m_P-1=2\) at all eight points; equivalently, all eight points on \(H\) are triple points. The Ziegler restriction to \(H\) is thus a multiarrangement of eight points, each of multiplicity two. By Lemma~\ref{lem:const2}, its exponents are \((8,8).\) Since \(d_1=8\), Lemma~\ref{lem:yoshinaga} then implies that \(\cA_2\) is free, contradicting the choice of \(\cA_2\). Hence every eight-point line contains a double point, proving \textup{(W4)}.
\end{proof}

\begin{lemma}\label{lem:const2}
A multiarrangement of $s\ge2$ distinct points on $\PP^1_\C$, each of multiplicity two, has exponents $(s,s)$.
\end{lemma}

\begin{proof}
This is \cite[Proposition~1.23(iii)]{Yoshinaga2014}.
\end{proof}

\begin{definition}\label{def:window}
A lattice satisfying (W1)--(W5) is called a \emph{window lattice}. A realisation is always understood to have exactly that lattice, with no extra incidences.
\end{definition}

The conditions \textup{(W1)--(W5)} are necessary for the lattice of a counterexample, but they do not characterize counterexamples. In particular, a
window lattice may admit only free realisations, only non-free realisations, or realisations with different freeness behaviour. The purpose of the conditions
above is therefore to delimit the finite combinatorial search; freeness must still be checked separately on the resulting realisation spaces.

We record two elementary identities that will be used repeatedly in this search. Let \(t_k\) denote the number of intersection points of multiplicity \(k\), let \(p\) be the total number of intersection points, and set
\[
X=\sum_P\binom{m_P-1}{2},
\qquad
W=\sum_P\binom{m_P-2}{2}.
\]
Double points contribute zero to both sums. Since
\[
\binom{m_P-1}{2}
=
\binom{m_P}{2}-(m_P-1),
\]
counting unordered pairs of lines gives \(X=\binom{n}{2}-b_2\). Moreover,
\[  
\binom{m_P-1}{2}-\binom{m_P-2}{2}=m_P-2,
\]
so that \(X-W =\sum_P(m_P-2) =b_2-p.\)

For a window lattice, \(n=17\) and \(b_2=80\). Hence
\(X=56\) and \(p=W+24\). Also,
\[
\sum_H s_H
=
\sum_P m_P
=
b_2+p,
\]
and therefore \(\sum_H(8-s_H)= 17\cdot8-(b_2+p)=32-W.\) Finally, since \textup{(W2)} gives \(m_P\le6\), we have
\[
X=t_3+3t_4+6t_5+10t_6,
\qquad
W=t_4+3t_5+6t_6,
\]
while \(b_2=t_2+2t_3+3t_4+4t_5+5t_6=80\). Eliminating \(t_3\) yields
\[
t_2=3W-32-t_5-3t_6.
\]
Thus the identities can be summarized as
\begin{equation}\label{eq:windowcounts}
X=56,\qquad
p=W+24,\qquad
\sum_H(8-s_H)=32-W,\qquad
t_2=3W-32-t_5-3t_6.
\end{equation}
In particular, \(t_2\ge0\) implies \(W\ge11\).

\section{Four deletions and a free thirteen-line core}\label{sec:RR}

Let $\cA$ be a non-free realisation of a window lattice $L$. By Theorem~\ref{thm:17}, it is nearly free with exponents $(7,10)$. We now compare it with the arrangement left after deleting four lines. This comparison has two parts: a combinatorial count determines the Chern classes, and the minimal syzygy controls the sections of the smaller bundle.

\subsection{The deletion count}
For a set $S$ of $k$ lines in $\cA$, call an intersection point \emph{full} if every arrangement line through it belongs to $S$, and let $f(S)$ count such points. Define
\begin{equation}\label{eq:phi}
\varphi(S)=\sum_{H\in S}(s_H-7)-f(S).
\end{equation}
\begin{lemma}\label{lem:count}
For every $S$, \(b_2(\cA\setminus S)=80-7k-\varphi(S).\) In particular, for $k=4$ we have $b_2(\cA \setminus S)=52-\varphi(S)$ and $\varphi(S)\le4$.
\end{lemma}
\begin{proof}
By \eqref{eq:comb}, every intersection point \(P\) contributes \(m_P-1\) to \(b_2(\cA)\). Suppose that \(a_P\) of the \(m_P\) lines through \(P\) are deleted. If \(a_P<m_P\), then at least one line through \(P\) remains. If at least two remain, the contribution of \(P\) to \(b_2(\cA\setminus S)\)
is \((m_P-a_P)-1\), so its contribution decreases by
\[
(m_P-1)-\bigl((m_P-a_P)-1\bigr)=a_P.
\]
If exactly one line remains, then \(P\) is no longer an intersection point of \(\cA\setminus S\), so its new contribution is zero.

On the other hand, if all the lines through \(P\) are deleted, then \(a_P=m_P\), and \(P\) disappears completely from the deleted arrangement. Its contribution therefore decreases from \(m_P-1\) to
zero. These are precisely the points counted by \(f(S)\). Hence the total
decrease of \(b_2\) is
\[
\sum_P a_P-f(S).
\]

Now \(\sum_P a_P\) counts the pairs \((P,H)\) such that
\(H\in S\) and \(P\in H\). Counting the same pairs line by line gives
\[
\sum_P a_P=\sum_{H\in S}s_H.
\]
Therefore
\[
b_2(\cA\setminus S)
=
80-\left(\sum_{H\in S}s_H-f(S)\right).
\]
Since
\[
\varphi(S)
=
\sum_{H\in S}(s_H-7)-f(S)
=
\left(\sum_{H\in S}s_H\right) -7k-f(S),
\]
we have
\[
\sum_{H\in S}s_H-f(S)=7k+\varphi(S),
\]
and consequently \(b_2(\cA\setminus S) = 80-7k-\varphi(S).\)

If \(k=4\), this becomes \(b_2(\cA\setminus S)=52-\varphi(S).\) Moreover, by \textup{(W3)}, every line satisfies \(s_H\leq8\), so
\[
\varphi(S)
=
\sum_{H\in S}(s_H-7)-f(S)
\leq4,
\]
because \(f(S)\geq0\).
\end{proof}

The correction $f(S)$ is essential: an intersection of two deleted lines alone contributes one, rather than two, to the drop in $b_2$. Equality $\varphi(S)=4$ means that all four lines have eight points and no intersection point is full. Following the computation, we call such an $S$ a \emph{$4$-chain}.

\subsection{Why the extremal core is free}
\begin{lemma}\label{lem:mdrS}
For every set $S$ of four lines, $\mdr(\cA\setminus S)\ge6$.
\end{lemma}
\begin{proof}
Let \(\theta\) be a minimal Jacobian syzygy of \(\cA\), of degree \(r = 7\), and let \(Z\) be its length-one zero scheme as a section of \(\cE_{\cA}(7)\). The associated section sequence is
\[
0\longrightarrow
\cO(-7)
\xrightarrow{\theta}
\cE_{\cA}
\longrightarrow
\cI_Z(-9)
\longrightarrow0.
\]
After twisting by \(k\), we have \(H^0\bigl(\cI_Z(k-9)\bigr)=0\) for \(k\leq9.\) It follows that every Jacobian syzygy of \(\cA\) of degree at most \(9\) is a polynomial multiple of \(\theta\); in particular, the
degree-seven syzygy is unique up to scalar.

Suppose now that \(\cA\setminus S\) admits a nonzero Jacobian syzygy of degree at most five. Multiplying it, if necessary, by a homogeneous form, we obtain a nonzero syzygy \(\psi\in D_0(\cA\setminus S)\) of degree five. Put
\[
Q_S=\prod_{H\in S}\alpha_H.
\]
Since \(|S|=4\), the derivation \(Q_S\psi\) has degree nine, and multiplication by the defining equations of the deleted lines makes it logarithmic along the whole arrangement \(\cA\). Using the decomposition
\[
D(\cA)=R\eu\oplus D_0(\cA)
\]
and the description \(D_0(\cA)_9 = R_2\theta\) obtained above, we may write
\[
Q_S\psi=g\eu+q\theta,
\qquad
\deg q=2.
\]

As in the foliation construction, let
\[
\Omega=dx\wedge dy\wedge dz,
\qquad
\omega_\eta:=\iota_{\eu}\iota_\eta\Omega
\]
for a homogeneous derivation \(\eta\); in particular, \(\omega=\omega_\theta=\iota_{\eu}\iota_\theta\Omega.\)
Applying this projectivization map to the preceding identity gives \(Q_S\omega_\psi=q\omega,\)
since \(\iota_{\eu}\iota_{g\eu}\Omega=0.\)

By Lemma~\ref{lem:sat}, the coefficients of \(\omega\) have no common factor. Let \(H\in S\). Since \(\alpha_H\mid Q_S\), the identity \(Q_S\omega_\psi=q\omega\) shows that \(\alpha_H\) divides every coefficient of \(q\omega\). Hence
\[
\alpha_H\mid q
\qquad\text{for every }H\in S.
\]
The four lines in \(S\) are distinct, so their defining linear forms
are pairwise non-associated. Consequently, \(Q_S\mid q.\) But \(\deg Q_S=4\) whereas \(\deg q=2\), thus \(q=0\).

It follows that \(Q_S\omega_\psi=0,\) and hence \(\omega_\psi=0\). By the kernel description of the
projectivization map used in the proof of Lemma~\ref{lem:sat},
\(\psi\) must be an Euler multiple \(\psi=h\eu.\)
This is impossible for a nonzero Jacobian syzygy. Indeed,
\(D_0(\cA\setminus S)\cap R\eu=0\). Thus
\(\cA\setminus S\) has no nonzero Jacobian syzygy of degree at most
five.
\end{proof}

\begin{theorem}[free cores]\label{thm:RR}
Let $\cA$ be a non-free realisation of a window lattice and let $S$ be a set of four lines from $\cA$. Then $\cA\setminus S$ is free with exponents $(6,6)$ if and only if $\varphi(S)=4$.
\end{theorem}
\begin{proof}
Write $\cE=\cE_{\cA\setminus S}$. The remaining arrangement has thirteen lines, so \eqref{eq:chern}, \eqref{eq:twist}, and Lemma~\ref{lem:count} give
\[
c_1(\cE(6))=0,\qquad c_2(\cE(6))=4-\varphi(S).
\]
Suppose $\varphi(S)=4$. Riemann--Roch gives $\chi(\cE(6))=2$, whereas duality and Lemma~\ref{lem:mdrS} give $h^2(\cE(6))=h^0(\cE(3))=0$. Thus $\cE(6)$ has a nonzero section $\nu$. Together with Lemma~\ref{lem:mdrS}, this gives $\mdr(\cA \setminus S) = 6$, so the section $\nu$ is of minimal degree and has no divisorial zero. Its zero scheme is finite and has length $c_2(\cE(6))=0$. This forces a split short exact sequence
\[
0\longrightarrow\cO \xlongrightarrow{\nu} \cE(6)\longrightarrow\cO\longrightarrow0,
\]
showing $\cE=\cO(-6)^2$. Conversely, freeness with exponents $(6,6)$ gives $b_2=36+12=48$, hence $\varphi(S)=4$, by Lemma~\ref{lem:count}.
\end{proof}

This also constrains the core lattice: every line of a free arrangement with exponents $(6,6)$ contains at most seven distinct intersection points \cite[Cor.~1.2(1)]{Abe14}. Moreover, freeness is combinatorial for thirteen lines \cite[Thm.~3.1]{DIM19}, so every realisation of this core lattice is free. These conclusions concern cores arising from a non-free realisation of a window lattice; the nearly free hypothesis was used in Lemma~\ref{lem:mdrS}.

\subsection{Recovering the deleted lines}
To reconstruct a seventeen-line arrangement from one of its thirteen-line
cores, we need to know that the deleted lines do not introduce any
continuous freedom once the remaining lines are fixed. The bounds
\textup{(W2)} and \textup{(W3)} imply precisely this: each deleted line
must pass through at least two intersection points already determined by
the lines that are present.

\begin{lemma}\label{lem:twomult}\label{lem:determined}
Let \(\cA\) be a seventeen-line arrangement satisfying \(m_P\leq6\) and $s_H\leq8$ for every intersection point \(P\) and every line \(H\in\cA\). Let
\(S\subseteq\cA\) contain at most four lines. Starting from
\(\cA\setminus S\), the lines of \(S\) can be restored in any order,
and at each step the next missing line is the join of two intersection
points of the arrangement already present.

Consequently, once the labelled incidence data of the full arrangement
are prescribed, a realisation of \(\cA\setminus S\) uniquely determines the missing lines. In particular, every line \(H\in\cA\) contains at least two intersection points \(P\) with \(m_P\geq3.\)
\end{lemma}

\begin{proof}
Consider a missing line \(H\) at some stage of the reconstruction.
Since at most four lines were deleted, at least thirteen other lines are
already present. All of them meet \(H\), and by the hypothesis
\(s_H\leq8\), these intersections occur at at most eight distinct
points of \(H\).

We claim that at least two of the intersection points of \(H\) with the lines already present are themselves intersections of at least two present lines. If there were none, the at least thirteen present lines would meet \(H\) in at least thirteen distinct points, contradicting \(s_H\leq8\). If there was exactly one, at most five present lines could pass through that point, since $m_P \leq 6$, while each of the remaining at most seven points could contain only one present line. Thus at most \(5+7=12\) present lines could meet \(H\), again a contradiction with the number of lines.

Thus \(H\) contains at least two points that are already intersection
points of the current arrangement. Their positions are determined by
the lines already present, and \(H\) is their unique join. This argument
can be repeated after each line is restored, so the lines of \(S\) may
be added back in any order.

For the final assertion, apply the same argument with \(S=\{H\}\). The two intersection points already present on \(H\) each lie on at
least two other arrangement lines; after \(H\) is included, their
multiplicity is therefore at least three.
\end{proof}

Here \emph{labelled} means that the correspondence between the lines in a
realisation and the line labels in the prescribed intersection lattice
is fixed. Thus the assertion does not identify realisations modulo
permutations of the lines; projectively equivalent labelled
realisations are still identified in the usual way. The conditions
requiring a realisation of a core to extend to the prescribed larger
arrangement are incidence conditions, given by vanishing and
non-vanishing of determinants, and are therefore defined over
\(\Q\).

The deficit behaves monotonically during this reconstruction.

\begin{lemma}\label{lem:monotone}
Fix a set \(S\) of lines. If the ambient arrangement is enlarged by
adding further lines while the lines of \(S\) are kept fixed, then
\(\varphi(S)\) cannot decrease.
\end{lemma}

\begin{proof}
Recall that
\[
\varphi(S)
=
\sum_{H\in S}(s_H-7)-f(S).
\]
Adding lines cannot decrease the number \(s_H\) of intersection points
on any fixed line \(H\in S\). Moreover, it cannot create a new full
point with respect to \(S\): any new intersection point involving an added line lies
on a line outside \(S\). An existing full point may cease to be full if
a new line passes through it, so \(f(S)\) can only decrease. Both
effects can only increase \(\varphi(S)\).
\end{proof}

The following arrangement illustrates how the deletion invariant,
Theorem~\ref{thm:RR}, and the reconstruction lemma interact in a
concrete example.

\begin{example}\label{ex:E17}\label{ex:E17chains}
Let $\zeta \in \mathbb{C}$ with \(\zeta^2-\zeta+1=0\) and consider the arrangement
\begin{align*}
E_{17}
={}\{z=0\}
 &\cup\{y=az:a\in\{0,1,\zeta-1,\zeta\}\}\cup\{x=bz:b\in\{0,1,\zeta,1-\zeta\}\}\\
&\cup\{x+y=cz:c\in\{0,1,\zeta,\zeta+1\}\}\cup\{x+\zeta y=dz:d\in\{0,1,\zeta-1,\zeta\}\}.
\end{align*}
Its intersection multiplicities are
\[
(t_2,t_3,t_4,t_5,t_6)=(24,8,8,4,0).
\]
Among its seventeen lines, one contains four distinct intersection
points, eight contain seven, and eight contain eight; we abbreviate
this line-point distribution by \(4\,7^8\,8^8.\)

The exact computation in
\cite{TC17} shows that \(E_{17}\)
is nearly free with exponents \((7,10)\) and that its lattice contains
thirty \(4\)-chains. By Theorem~\ref{thm:RR}, deleting any one of these
\(4\)-chains produces a free thirteen-line core with exponents
\((6,6)\).

The arrangement \(E_{17}\) is obtained from the reflection arrangement
of \(G_{26}\) by deleting four lines and appears as the first lattice
in the final list. Its lattice has only a conjugate pair of
realisations, both with the same nearly free behaviour. Thus
\(E_{17}\) is not a counterexample to Terao's conjecture; rather, it
provides a concrete example of the reduction from a non-free
seventeen-line arrangement to free thirteen-line cores.
\end{example}

\section{Finding a useful deletion}\label{sec:cover}

The free-core theorem identifies the best possible four-line deletion. For the reduction we need less: a deletion with positive deficit. Its existence is the only combinatorial assertion in the reduction that we prove by a finite computation.

\begin{lemma}[covering lemma; computer-assisted]\label{lem:cover}
If a lattice satisfying \textup{(W1)--(W4)} has an eight-point line, some set $S$ of four lines satisfies $\varphi(S)\ge1$.
\end{lemma}

\begin{proof}
Fix an eight-point line $H_0$. Its eight points partition the other sixteen lines into pencils $G_1,\ldots,G_8$, of sizes $k_i=m_{P_i}-1$. Thus $1\le k_i\le5$, $\sum k_i=16$, and some $k_i=1$. There are fourteen such unordered vectors. Put $\delta_H=8-s_H$ and $\Delta=\sum_H\delta_H=32-W$; in particular $\delta_{H_0}=0$ and $\Delta\le21$.

Suppose that every four-line set has non-positive deficit. For lines $A,B$ in different pencils let $e_{AB}$ indicate that their intersection is double. For lines $A,B,C$ in three different pencils let $\tau_{ABC}$ indicate that their intersection is triple. Set all other such indicators to zero. The full points of $\{H_0,A,B,C\}$ are precisely the pencil centres whose whole pencil is selected, the indicated double points, and the indicated triple point. Therefore
\begin{equation}\label{eq:triple}
\delta_A+\delta_B+\delta_C+
\sum_i[G_i\subseteq\{A,B,C\}]+e_{AB}+e_{AC}+e_{BC}+\tau_{ABC}\ge4.
\end{equation}
Here brackets denote an indicator. A pair can lie in at most one indicated double or triple point, so $e_{AB}+\sum_C\tau_{ABC}\le1$; it can instead belong to a point of higher multiplicity.

Write $b'_X=\sum_Ye_{XY}$ for the number of double points of $X\in G_i$ off $H_0$, and $s=8-\delta_X$. The $s-1$ intersection points off $H_0$ account for $16-k_i$ other lines. A double point contributes one and every other point between two and five, giving
\begin{equation}\label{eq:coverbounds}
2s-18+k_i\le b'_X\le
U_X:=\min\left\{\left\lfloor\frac{5s-21+k_i}{4}\right\rfloor,
s-3+[k_i\ge2]\right\}.
\end{equation}
The second upper bound uses the two multiple points on $X$ supplied by Lemma~\ref{lem:twomult}. Also $b'_X\ge1$ when $s=8$ and $k_i\ge2$, by (W4). Finally, with $u=\#\{i:k_i=1\}$ and $t_j(H_0)=\#\{i:k_i=j-1\}$,
\begin{equation}\label{eq:coverglobal}
\begin{split}
\Delta&\le32-\max\left\{11,\sum_i\binom{k_i-1}{2}\right\},\\
u+\sum_{\{X,Y\}}e_{XY}&\le64-3\Delta-t_5(H_0)-3t_6(H_0).
\end{split}
\end{equation}
The pair sum is unordered. These inequalities follow from $W\ge11$ and $t_2=64-3\Delta-t_5-3t_6$.

The contradiction is certified in two finite steps, recorded in \path{work/n17/cover/ilp/} of \cite{TC17}. First enumerate the assignments $0\le\delta_X\le4$, modulo permutations within pencils and between equal-sized pencils. Retain the consequences of \eqref{eq:coverbounds}--\eqref{eq:coverglobal} involving only the deficits, together with the following consequence of \eqref{eq:triple}: for $A,B$ in the same pencil of size $k\ge2$, summing over the fourteen remaining choices of $C$ gives
\[
U_A+U_B\ge56-\Delta-13(\delta_A+\delta_B)-u-14[k=2]-[k=3].
\]
It suffices for pruning to impose this on the two lines of least deficit in each pencil of size at least three and on the unique pair when $k=2$. The program \code{level1.py} leaves $2{,}471$ assignments.

For each assignment, relax the remaining indicators to $0\le x_j\le1$ and write the inequalities as $\sum_jc_{ij}x_j\ge b_i$. The file of certificates supplies rational multipliers $y_i\ge0$ with
\begin{equation}\label{eq:farkas}
\sum_j\max\left\{0,\sum_i y_ic_{ij}\right\}<\sum_i y_i b_i.
\end{equation}
Multiplying and adding the inequalities contradicts \eqref{eq:farkas}, because each $x_j$ lies between zero and one. GLPK \cite{glpk} proposes the multipliers; \code{farkas.py} verifies the strict inequality with exact rational arithmetic in all $2{,}471$ cases. Thus the supposed lattice cannot exist.
\end{proof}

Put $\varphi_{\max}(L)=\max_{|S|=4}\varphi(S)$, and call the deletion of $S$ from $L$ a \emph{core} when $S$ attains this maximum.

\begin{corollary}\label{cor:cores}
A window lattice with an eight-point line has $1\le\varphi_{\max}\le4$. Every core has thirteen lines, $b_2=52-\varphi_{\max}\in\{48,49,50,51\}$, multiplicities at most six, and three to eight points on each line. No two points of multiplicities $\{6,5\}$ or $\{6,6\}$ are unjoined. If the lattice has a non-free realisation and $\varphi_{\max}=4$, every core is free with exponents $(6,6)$ and has at most seven points on each line.
\end{corollary}

\begin{proof}
The deficit bounds and the value of $b_2$ follow from Lemmas~\ref{lem:cover} and~\ref{lem:count}. Deletion preserves the multiplicity bound and cannot increase a line's point count. Twelve other core lines meet a given line in groups of at most five, so at least three points remain. A multiplicity-six point of the core cannot acquire a deleted line; hence an unjoined pair containing such a point was already unjoined before deletion, and (W5) excludes the stated pairs. When a non-free realisation exists, Theorem~\ref{thm:17} makes it nearly free, and Theorem~\ref{thm:RR} makes a maximal-deficit core free. Abe's bound \cite[Cor.~1.2(1)]{Abe14} then gives at most seven points on each core line.
\end{proof}

\section{The finite verification}\label{sec:computation}

The reduction has replaced a large search on seventeen lines by two manageable tasks: eliminate lattices without an eight-point line, and extend thirteen-line cores. This section specifies the finite searches, explains why they cover every possible counterexample, and records their outputs. Code and exact records are available in \cite{TC17}; the computational statements below refer to those records.

\subsection{Exact realisation spaces}\label{ss:filters}
We first describe the algebraic procedure used throughout. In the dual plane an arrangement is a configuration of labelled points $P_i$; its multiple points become prescribed collinear subsets, or \emph{blocks}. Choose four elements with no three in a block and normalise their points to
\[
(1:0:0),\quad(0:1:0),\quad(0:0:1),\quad(1:1:1).
\]
Every realisation admits this unique frame normalisation. Place each remaining point using as many prescribed collinearities with already placed points as possible: two distinct prescribed lines determine an intersection, one prescribed line gives $P_a+uP_b$, and none gives $(u:v:1)$. These choices do not lose a realisation. In the one-parameter case the missing endpoint $P_b$ would coincide with an existing point. In the two-parameter case a point on $z=0$ would be collinear with the first two frame points, so it would either have been placed with fewer parameters or violate the prescribed lattice.

Let $I$ be the ideal of prescribed determinant equations and let $D$ be the product of the determinants required to be non-zero. Thus the frame-fixed realisation space is
\begin{equation}\label{eq:realspace}
R(L)=V(I)\setminus V(D).
\end{equation}
The non-vanishing conditions exclude both extra collinearities and coincident points. All coefficients are rational. A minimal prime $\mathfrak p$ of $I$ contributes a component precisely when no factor of $D$ lies in $\mathfrak p$. This is a statement about components over $\Q$; a zero-dimensional component may contain several complex points, forming one Galois orbit.

Before computing minimal primes, inexpensive exact tests reject many empty spaces. An identically zero required determinant or a non-zero constant incidence equation is decisive. With one parameter, take the gcd of the equations and remove every factor shared with a degeneracy polynomial; a constant remainder means that no non-degenerate solution exists. With two parameters $u,v$, every solution satisfies the equations independent of $v$ and the $v$-resultants of pairs of equations that both involve $v$. A constant gcd of these univariate conditions excludes a solution. Otherwise its irreducible factors $f(u)$ can be tested separately over $\Q[u]/(f)$ by the one-variable procedure in $v$. Common curve factors of the equations are treated separately and removed only when a degeneracy polynomial vanishes identically on them. An inconclusive test always defers the lattice.

These filters use exact polynomial arithmetic in FLINT \cite{flint}. In particular, two equations independent of $v$ are retained as equations in $u$, not fed to the resultant operation: their resultant with respect to an absent variable would give an invalid rejection. Macaulay2 \cite{M2} computes the minimal primes for every deferred lattice. For a zero-dimensional component it also supplies coordinates $u_i=g_i(T)$ with $f(T)=0$, giving an exact realisation over $K=\Q[T]/(f)$.

Freeness is a separate linear-algebra calculation. On a chosen arrangement line, a degree-$k$ multirestriction derivation is a pair of binary forms $(a,b)$ satisfying $\beta_P^{m_P-1}\mid\theta(\beta_P)$ at every point. Equating remainders to zero gives a matrix with $2(k+1)$ columns. For $k=d_1-1$, Lemma~\ref{lem:yoshinaga} says that freeness is equivalent to full column rank. For the final seventeen-line arrangements we compute this rank over each component's number field, thereby deciding all conjugate realisations simultaneously. No inference about special members of a family is made from its generic rank.

\subsection{The branch without an eight-point line}\label{ss:no8}
The identities in \eqref{eq:windowcounts} give
\begin{equation}\label{eq:no8budget}
\sum_H(7-s_H)=15-W.
\end{equation}
When every line has at most seven points, the left side is non-negative. A line with four points uses three units of this budget, and a line with five points uses two.

\begin{lemma}\label{lem:no8}
Let $L$ be a window lattice with no eight-point line. If $L$ has a four-point line, its four points all have multiplicity five, every other line has seven points, and the remaining multiple points give a triangle decomposition of $K_{4,4,4,4}$. If $L$ has a five-point line, then $W\le13$.
\end{lemma}

\begin{proof}
For a four-point line $H_0$, \eqref{eq:no8budget} gives $W\le12$. Writing its multiplicities as $k_i+1$, we have $\sum_{i=1}^4k_i=16$, $k_i\le5$, and $W\ge\sum_i\binom{k_i-1}{2}$. Among the five possible unordered vectors, only $(4,4,4,4)$ has this last sum at most twelve; the others give $13,14,15,18$. Thus $W=12$, every other line has seven points, and no further point has multiplicity at least four. The counting identities give $t_2=0$ and $t_3=32$. The four pencils off $H_0$ form the four parts of $K_{4,4,4,4}$: each triple point covers exactly three edges, and all $96$ edges are covered once. The assertion for a five-point line follows directly from \eqref{eq:no8budget}.
\end{proof}

\begin{theorem}[computer-assisted]\label{thm:no8}
No window lattice without an eight-point line is realisable over $\C$.
\end{theorem}

\begin{proof}
There are three exhaustive cases.

\emph{A four-point line.} An exact-cover enumeration, with isomorphism decided by \texttt{nauty} \cite{MP14}, gives $23$ triangle decompositions of $K_{4,4,4,4}$. The filters reject nine; Macaulay2 rejects the other fourteen.

\emph{A five-point line and no four-point line.} Write the five pencil sizes as $(k_1,\ldots,k_5)$, where $\sum k_i=16$ and $k_i\le5$. Of the eleven possible vectors, $(5,5,4,1,1)$ has $W(H_0)=15>13$. Enumerate the points of multiplicity at least four off $H_0$, respecting the budget $W-W(H_0)\le13-W(H_0)$, and choose the point counts $5\le s_X\le7$. If $X\in G_i$ contains $h$ such points of multiplicities $m_1,\ldots,m_h$, its remaining triple and double point counts are forced:
\[
\mathrm{tr}_X=17-k_i-s_X-h-\sum_j(m_j-3),\qquad
\mathrm{db}_X=s_X-1-h-\mathrm{tr}_X.
\]
Both must be non-negative. Completing the remaining pairs by triples subject to these quotas is an exact-cover search. The five vectors $(5,5,3,2,1)$, $(5,5,2,2,2)$, $(5,4,4,2,1)$, $(5,4,3,3,1)$ and $(5,4,3,2,2)$ admit no quota system. Table~\ref{tab:five} records the other five. Every resulting lattice is excluded by the filters or by Macaulay2. The enumeration fixes $H_0$ when identifying partial configurations; forgetting this distinguished line would not justify completeness. The programs are \code{pencil5c.py} and \code{p5excl.py} under \path{work/n17/n8zero/}.

\emph{Only six- and seven-point lines.} Generate the configurations of blocks of size at least four, complete by triples, and identify isomorphic lattices by canonical incidence graphs. The conditions $b_2=80$, $m\le6$, $6\le s_H\le7$ and (W5) leave $1{,}958{,}407$ lattices. The corrected filters reject $1{,}943{,}787$; Macaulay2 rejects the remaining $14{,}620$, with unit incidence ideal in $1{,}301$ cases and only degenerate components in $13{,}319$. The programs are \path{work/tc15/bigblocks2.py} and \code{pipeline.py}.
\end{proof}

\begin{table}[ht]
\centering\small
\begin{tabular}{lrrr}
\toprule
Pencil sizes & Heavy configurations & Lattices & Deferred to Macaulay2\\
\midrule
$(5,3,3,3,2)$ & 2,303 & 0 & 0\\
$(4,4,4,3,1)$ & 87 & 122 & \multirow{2}{*}{71}\\
$(4,4,4,2,2)$ & 863 & 232 & \\
$(4,4,3,3,2)$ & 49,565 & 3,124 & 464\\
$(4,3,3,3,3)$ & 192,133 & 7,723 & 355\\
\bottomrule
\end{tabular}
\caption{Five-point-line enumeration. None of the lattices is realisable.}\label{tab:five}
\end{table}

\subsection{Enumerating the cores}\label{ss:cores}
We now turn to lattices with an eight-point line. A lattice on $q$ lines is encoded by blocks of size at least three, with no pair in two blocks; an uncovered pair is a double point. Thus
\[
X=\sum_B\binom{|B|-1}{2}=\binom q2-b_2,
\qquad w(x)=\sum_{B\ni x}(|B|-2)=q-1-s_x.
\]
These identities translate the core bounds into integer constraints on blocks.

\begin{proposition}[computer-assisted]\label{prop:cores}
The thirteen-line classes specified in Corollary~\ref{cor:cores}, with $3\le s_H\le7$ for $b_2=48$ and $3\le s_H\le8$ for $b_2=49,50,51$, have the counts in Table~\ref{tab:cores}. In all four classes every multiplicity is at most six and unjoined pairs of multiplicities $\{6,5\}$ and $\{6,6\}$ are excluded.
\end{proposition}

\begin{proof}
For each $b_2$, prescribe $X=78-b_2$ and $4\le w(x)\le9$, with the lower bound raised to five when $b_2=48$. First generate the configurations of blocks of size at least four by canonical augmentation \cite{McK98}. The adapted program \code{heavyaug} of Anghel and Stan \cite{AS26} gives respectively $258,274,284,285$ configurations. Their program \code{lsgen3} completes these by triple blocks; uncovered pairs remain double points. It produces respectively $213{,}339$, $1{,}820{,}363$, $12{,}276{,}119$ and $65{,}691{,}853$ completions. Canonical labelling of the incidence graph, with line-vertices and point-vertices distinguished, removes repetitions and gives the four counts in the table. The exact procedure of Section~\ref{ss:filters} then gives the realisation counts. All $186{,}771$ deferred core lattices are resolved by Macaulay2.
\end{proof}

\begin{table}[ht]
\centering\small
\begin{tabular}{crrrrrr}
\toprule
$b_2$ & Lattices & Deferred & Realisable & Rigid & Curves & Dim.\ $\ge2$\\
\midrule
48 & 21,164 & 195 & 85 & 80 & 5 & 0\\
49 & 154,172 & 2,950 & 1,147 & 1,050 & 83 & 14\\
50 & 738,023 & 26,093 & 9,162 & 7,985 & 1,091 & 86\\
51 & 2,916,332 & 157,533 & 58,871 & 49,118 & 9,001 & 752\\
\midrule
Total & 3,829,691 & 186,771 & 69,265 & 58,233 & 10,180 & 852\\
\bottomrule
\end{tabular}
\caption{The core classes. The last three columns record the dimension categories of their realisation spaces.}\label{tab:cores}
\end{table}

\subsection{Extending every realisation}\label{ss:extension}
The extension step must cover whole realisation spaces, not merely one convenient representative of each core lattice. Lemma~\ref{lem:determined} supplies the finite geometric operation: at each step, add a line through two existing intersection points. Every desired seventeen-line realisation is obtained after four such additions. The implementation separates three cases according to the dimension of the current component.

\emph{A zero-dimensional component.} Work over its number field $K$. List all joins of distinct intersection points, compute their incidences, and recurse on every admissible new line. Arithmetic is exact. Reduction modulo a prime may certify a value to be non-zero; a value whose reduction vanishes is checked in $K$. This shortcut never certifies equality merely from a modular coincidence.

\emph{A curve component.} Let $C$ be an irreducible curve component of the original non-degenerate incidence space. A plane model $F(T,Y)=0$ represents its function field, and the coordinates of the arrangement are rational functions on that model. To ensure completeness, certify dense open subsets $U\subseteq C$ and $V\subseteq V(F)$ on which the birational maps are inverse isomorphisms. Compute the finite complement $C\setminus U$ in the original incidence coordinates, retaining all non-vanishing conditions that define the prescribed lattice, and search every point of this complement over its number field. This includes valid realisations omitted by the affine plane model, including those mapped to infinity. Run the generic search over the function field on $U$. An incidence test is either identically satisfied or can change at a zero of a resultant with $F$. Such exceptional values must also be searched: Macaulay2 finds the corresponding points of the original incidence variety, and each is passed to the exact number-field search. Zeros of denominators are included as well. The implementation restricts these tests to specialisations that can make a candidate line admissible: extra coincidences can reduce its point count, and all exceptional roots for an already admissible line are retained. Thus the generic calculation and the exceptional-point calculations are distinct obligations; the latter cannot be replaced by sampling the parameter.

\emph{Dimension at least two.} Generate possible new-line incidence patterns combinatorially. The old multiple points chosen on the new line must be pairwise unjoined by an old line; their multiplicities are at most five. All other old lines meet the new line in separate new double points. Retain only patterns compatible with the final bounds $b_2=80$, $m\le6$ and $s_H\le8$. Screen and decompose each resulting lattice exactly, then return its components to the appropriate case above. Iterating reaches seventeen lines. This also provides a fallback for an entire curve component when its chart coverage cannot be certified or its function-field calculation fails or exceeds its time limit; a failure is never interpreted as non-realisability.

One further reduction avoids repeating work between core classes. In a search starting from $B$, discard a partial extension if it already contains a four-line set with deficit larger than $52-b_2(B)$. To see why this is sufficient for our targets, choose an enumerated core $B$ obtained from $L$ by deleting $S$, with $\varphi(S)=\varphi_{\max}(L)$. Every partial arrangement in its reconstruction has all four-line deficits at most $\varphi_{\max}(L)$ by Lemma~\ref{lem:monotone}. This reconstruction is therefore never discarded. The search processes the four classes in increasing order of $b_2$.

The programs \code{extend.py}, \code{famext.py}, \code{famspecial2.py}, \code{planemodel.py} and \code{dim2gen.py} under \path{work/n17/comp/} perform the extension searches described here, with completeness requiring the chart coverage described above. In the initial run, $297$ of $14{,}983$ curve searches and two timed-out child computations were routed to the combinatorial fallback; $146{,}345$ deferred children with fourteen to sixteen lines were resolved by Macaulay2. The orchestration records every component and its outcome. The recorded final check verifies completion of the fallbacks and processing of the exceptional points generated by these searches.

\subsection{The output and the proof of the theorem}\label{ss:result}
The resulting list is small enough for a second exact calculation from its incidence data, independent of the route by which its members were found.

\begin{theorem}[computer-assisted]\label{thm:sixteen}
Let $L$ be a complex-realisable window lattice with a core in the classes of Proposition~\ref{prop:cores}. Then $L$ is one of the sixteen lattices in Table~\ref{tab:sixteen}. Each has one zero-dimensional component over $\Q$: a single projective realisation over $\Q$, or a conjugate pair over $\Q(\sqrt{-3})$. Every realisation is nearly free with exponents $(7,10)$. Each lattice has an eight-point line and a four-line set of deficit four.
\end{theorem}

\begin{proof}
Theorem~\ref{thm:no8} supplies an eight-point line, and by hypothesis a maximizing core belongs to one of the enumerated classes. Lemma~\ref{lem:determined} reconstructs every realisation from that core. Its reconstruction is never removed by deficit pruning, as just shown. After identification of isomorphic outputs, exactly sixteen lattices remain. Their complete incidence data and exact verdicts are in \path{work/n17/run/checked_L.jsonl}. For each, the direct incidence-ideal calculation gives one component of degree one or two, and the multirestriction matrix in degree seven has rank fifteen out of sixteen over its number field. Thus no realisation is free by Lemma~\ref{lem:yoshinaga}; Theorem~\ref{thm:17} makes each nearly free. The line counts and deficit-four sets are checked directly from the output lattices.
\end{proof}

\begin{table}[ht]
\centering\small
\begin{tabular}{lcc@{\qquad}lcc}
\toprule
Lattice & $(t_2,t_3,t_4,t_5,t_6)$ & $s_H$ & Lattice & $(t_2,t_3,t_4,t_5,t_6)$ & $s_H$\\
\midrule
$\Gamma_1$ & $(24,8,8,4,0)$ & $4\,7^8\,8^8$ & $\Gamma_9$ & $(25,11,3,6,0)$ & $7^{11}\,8^6$\\
$\Gamma_2$ & $(23,5,13,2,0)$ & $5\,7^{10}\,8^6$ & $\Gamma_{10}$ & $(27,5,9,4,0)$ & $6^2\,7^7\,8^8$\\
$\Gamma_3$ & $(24,8,8,4,0)$ & $5\,7^9\,8^7$ & $\Gamma_{11}$ & $(26,8,6,5,0)$ & $6^2\,7^7\,8^8$\\
$\Gamma_4$ & $(23,5,13,2,0)$ & $6^3\,7^7\,8^7$ & $\Gamma_{12}$ & $(27,5,9,4,0)$ & $6^3\,7^5\,8^9$\\
$\Gamma_5$ & $(22,8,10,3,0)$ & $6^3\,7^7\,8^7$ & $R_{17}^{a}$ & $(23,10,8,2,1)$ & $5\,6^3\,7^3\,8^{10}$\\
$\Gamma_6$ & $(24,8,8,4,0)$ & $6\,7^{10}\,8^6$ & $R_{17}^{b}$ & $(20,18,2,2,2)$ & $5\,6^3\,7^3\,8^{10}$\\
$\Gamma_7$ & $(25,5,11,3,0)$ & $6^2\,7^8\,8^7$ & $N_1$ & $(23,15,3,2,2)$ & $6^4\,7^3\,8^{10}$\\
$\Gamma_8$ & $(24,8,8,4,0)$ & $6^2\,7^8\,8^7$ & $N_2$ & $(26,8,6,5,0)$ & $6^4\,7^3\,8^{10}$\\
\bottomrule
\end{tabular}
\caption{The sixteen output lattices; exponents in the $s_H$ column denote frequencies. The twelve $\Gamma_i$ are four-line deletions of the $G_{26}$ reflection arrangement and are defined over $\Q(\sqrt{-3})$. The other four are rational; $N_1,N_2$ were first found by the extension search.}\label{tab:sixteen}
\end{table}

\begin{proof}[Proof of Theorem~\ref{thm:main}]
Suppose a seventeen-line lattice supports a free arrangement and a non-free one. Theorem~\ref{thm:17} makes it a window lattice. It has an eight-point line by Theorem~\ref{thm:no8}, and its non-free realisation ensures that its cores satisfy Corollary~\ref{cor:cores}, including the seven-point bound when $b_2=48$. Theorem~\ref{thm:sixteen} then says that every realisation is non-free, a contradiction. Combining this with \cite{BK23,AS26} proves the assertion for at most seventeen lines. Central arrangements of rank at most two are free, giving the equivalent statement for planes in $\C^3$.
\end{proof}

\begin{remark}\label{rem:freelattices}
This proves the required assertion without classifying every realisable window lattice. A window lattice all of whose realisations are free could have $\varphi_{\max}=4$ and an eight-point line in every thirteen-line core; those $b_2=48$ cores were not enumerated. Such a lattice cannot be a counterexample. The classification above includes every window lattice having a non-free realisation, which is the scope needed for the proof.
\end{remark}

\appendix
\section{Records and verification}\label{sec:data}
The computation was completed on 21 October 2026 on a ten-core workstation and audited on 21--23 October. The repository \cite{TC17} contains the programs, certificates and verdicts; \path{work/n17/MANIFEST.md} records commands and SHA-256 hashes for 9.3~GB of regenerable data omitted from the repository. The directories \code{cover/}, \code{n8zero/}, \code{heavy/}, \code{comp/}, \code{results/}, \code{run/} and \code{audit/} contain respectively the covering certificates, no-eight-point searches, core enumeration, extensions, exact algebraic records, final outputs and audit. The shared filters remain under \path{work/tc15/}.

\subsection*{Checks and their scope}\label{ss:audit}
An independently implemented enumerator from \cite{TeraoHunter}, run without the unjoined-pair cut, reproduced the four core lists and added exactly the seven $b_2=51$ lattices excluded by that cut. Anghel and Stan's exact screening code \cite{AS26} was run on all $3{,}829{,}691$ cores and $541{,}254$ children; it agreed wherever it decided a lattice, and Macaulay2 checked those rejected only by our filters. The rigid extension was rerun with exact memoisation. A separate numerical implementation, with $50$-digit rechecks near zero, found identical lattice sets on $99{,}449$ rigid components; this is a cross-check, not the arithmetic basis of the proof. All $14{,}983$ curve components in the audit were resolved again using fresh plane models or combinatorial extension, with agreement on the $205$ components treated both ways. The covering and no-eight-point branches were rerun, and the sixteen final lattices were independently certified non-free with the code of \cite{AS26}.

The audit corrected four implementation defects without changing a verdict: forgetting the distinguished line in the five-point enumeration; asserting some pencil exclusions before running their checks; memoisation based on modular images without exact confirmation; and failure to record worker-process errors. The corrected filters, including the parameter-free resultant case described above, are used throughout the final records. No completed output was rejected by (W4), (W5), or the additional external-line test recorded in the source computation.

The remaining shared dependencies are explicit. The function-field extension with its exceptional specialisations, the combinatorial generation of child lattices, and Macaulay2's minimal-prime computations have not been independently reimplemented. The Farkas certificates certify the inequalities supplied to them; the derivation of that model is the mathematical argument in Section~\ref{sec:cover}. The new theoretical arguments have been checked by the authors. Software versions are Python 3.12, python-flint 0.9 \cite{flint}, pynauty 2.8, Macaulay2 1.26 \cite{M2}, nauty 2.9 \cite{MP14}, and GLPK 5.0 \cite{glpk}; the C enumerators were compiled with Apple clang 17.

\bibliographystyle{plain}
\bibliography{refs}
\end{document}